\UseRawInputEncoding
\documentclass[11pt]{article}

\DeclareMathAlphabet{\mathpzc}{OT1}{pzc}{m}{it}

\usepackage{amsmath}
\usepackage{amssymb}
\usepackage{amsfonts}
\usepackage{amsthm}
\usepackage{mathrsfs}
\usepackage{ifsym}
\usepackage{fontenc}
\usepackage{textcomp}
\usepackage{tipa}
\usepackage{enumerate}
\usepackage[pdftex]{color, graphicx}

\begin{document}

\title{{\bf The Hasse invariant of the Tate normal form $E_5$ and the class number of $\mathbb{Q}(\sqrt{-5l})$}}         
\author{Patrick Morton}        
\date{Jan. 5, 2021}          
\maketitle

\begin{abstract}  It is shown that the number of irreducible quartic factors of the form $g(x) = x^4+ax^3+(11a+2)x^2-ax+1$ which divide the Hasse invariant of the Tate normal form $E_5$ in characteristic $l$ is a simple linear function of the class number $h(-5l)$ of the field $\mathbb{Q}(\sqrt{-5l})$, when $l \equiv 2,3$ modulo $5$.  A similar result holds for irreducible quadratic factors of $g(x)$, when $l \equiv 1, 4$ modulo $5$. This implies a formula for the number of linear factors over $\mathbb{F}_p$ of the supersingular polynomial $ss_p^{(5*)}(x)$ corresponding to the Fricke group $\Gamma_0^*(5)$.
\end{abstract}

\section{Introduction}

Let $E_5(b)$ be the Tate normal form of an elliptic curve with a point of order $5$:
$$E_5(b): \ Y^2+(1+b)XY+bY=X^3+bX^2;$$
and let $\hat H_{5,l}(b)$ denote its Hasse invariant in characteristic $l$, as in \cite[p. 260]{mor} and \cite[p. 235]{mor8}:
\begin{align*}
\hat H_{5,l}(b) & = (b^4+12b^3+14b^2-12b+1)^r(b^2+1)^s(b^4+18b^3+74b^2-18b+1)^s\\
& \times b^{5n_l}(1-11b-b^2)^{n_l}J_l\left(\frac{(b^4+12b^3+14b^2-12b+1)^3}{b^5(1-11b-b^2)}\right).
\end{align*}
The roots of $\hat H_{5,l}(x)$ in characteristic $l \neq 5$ are of course the values of $b \in \overline{\mathbb{F}}_l$ for which $E_5(b)$ is supersingular.  Here, $n_l=[l/12], r = \frac{1}{2}(1-(-3/l)), s=\frac{1}{2}(1-(-4/l))$; and  $J_l(t)$ is the polynomial
$$J_l(t)\equiv \sum_{k=0}^{n_l}{{2n_l + s\atopwithdelims ( ) 2k + s}{2n_l-2k \atopwithdelims ( ) n_l-k}(-432)^{n_l-k}(t-1728)^k} \ \ (\textrm{mod} \ l).$$
Let $\varepsilon=\frac{-1+\sqrt{5}}{2}$ and $\bar \varepsilon = \frac{-1-\sqrt{5}}{2}$, and let $h(-5l)$ denote the class number of the quadratic field $K=\mathbb{Q}(\sqrt{-5l})$, where $l > 5$ is a prime number.  This paper is devoted to proving the following theorem, which was stated in \cite{mor8}.  \medskip

\noindent {\bf Theorem 1.1.} {\it Let $l>5$ be a prime.} \smallskip

\noindent A) {\it Assume $l \equiv 2,3$ mod $5$.  The number of irreducible quartics of the form $g(x) = x^4+ax^3+(11a+2)x^2-ax+1$ dividing $\hat H_{5,l}(x)$ over $\mathbb{F}_l$ is: }
\begin{enumerate}[i)]
\item {\it $\frac{1}{4}h(-5l)$, if $l \equiv 1$ mod $4$;}
\item {\it $\frac{1}{2}h(-5l)-1$, if $l \equiv 3$ mod $8$;}
\item {\it $h(-5l)-1$, if $l \equiv 7$ mod $8$.}
\end{enumerate}
B) {\it If $l \equiv 4$ mod $5$, the number of irreducible quadratic factors $k(x)=x^2+rx+s$ dividing $\hat H_{5,l}(x)$ over $\mathbb{F}_l$, with $r = \varepsilon ^5 (s-1)$ or $r = \bar \varepsilon^5 (s-1)$ in $\mathbb{F}_l$ is:}
\begin{enumerate}[i)]
\item {\it $\frac{1}{2}h(-5l)$, if $l \equiv 1$ mod $4$;}
\item {\it $h(-5l)-3$, if $l \equiv 3$ mod $8$;}
\item {\it $2h(-5l)-3$, if $l \equiv 7$ mod $8$.}
\end{enumerate}
C) {\it If $l \equiv 1$ mod $5$, the number of irreducible quadratic factors $k(x)=x^2+rx+s$ dividing $\hat H_{5,l}(x)$ over $\mathbb{F}_l$, with $r = \varepsilon ^5 (s-1)$ or $r = \bar \varepsilon^5 (s-1)$ in $\mathbb{F}_l$ is:}
\begin{enumerate}[i)]
\item {\it $\frac{1}{2}h(-5l)$, if $l \equiv 1$ mod $4$;}
\item {\it $h(-5l)-1$, if $l \equiv 3$ mod $8$;}
\item {\it $2h(-5l)-1$, if $l \equiv 7$ mod $8$.}
\end{enumerate}
\medskip

This theorem is a supplement to \cite[Thm. 6.2]{mor}, which gives the possible degrees of irreducible factors of $\hat H_{5,l}(x)$ (mod $l$), depending on the congruence class of $l$ (mod $5$).  In Part A, the only possible irreducible factors of $\hat H_{5,l}(x)$ are either $x^2+1$ or quartic factors $f(x)$ satisfying $x^4f(-1/x)=f(x)$; while in Part B, the possible factors are linear and quadratic; and in Part C, all irreducible factors are quadratic.  The class number $h(-5l)$ then shows up in counting the special factors of the form $g(x)$ or $k(x)$ in Theorem 1.1.  \medskip

In the paper \cite[Thms. 1.1 and 1.3]{mor8} I showed that the number of linear factors of $\hat H_{5,l}(x)$ (mod $l$) is a function of the class number $h(-l)$, for $l \equiv 4$ (mod $5$), and that the same holds for the number of quadratic factors of the form $x^2+ax-1$, when $l \equiv 1$ (mod $5$).  Thus, both class numbers $h(-l)$ and $h(-5l)$ are encoded in the factorization of $\hat H_{5,l}(x)$ over $\mathbb{F}_l$, at least when $l \equiv \pm 1$ (mod $5$).  I conjecture that the same is true for the Hasse invariant $\hat H_{p,l}(x)$ of the Tate normal form $E_p$ for a point of order $p$, i.e., that its factorization encodes the class numbers $h(-l), h(-pl)$, where $l \equiv \pm 1$ (mod $p$) and $p$ is any prime for which the coefficients of $E_p$ depend on a single parameter.  \medskip

Theorem 1.1 is an analogue for the Tate normal form $E_5$ of similar results for the Legendre normal form (or equivalently, the Tate normal form $E_4$ for a point of order $4$, see \cite[p. 255]{mor}) and the Deuring normal form.  The main theorems of \cite{mor4} and \cite{mor5} (and their proofs) show that the counts of certain irreducible quadratic factors of the Hasse invariants of these normal forms are functions of the class numbers $h(-2l), h(-3l)$ of the imaginary quadratic fields $\mathbb{Q}(\sqrt{-2l})$ and $\mathbb{Q}(\sqrt{-3l})$, respectively. \medskip

The proof of Theorem 1.1 generally proceeds along similar lines as the proof of the corresponding result for the Deuring normal form in \cite{mor5} and \cite{mor7}, and is a combination of theoretical arguments and computational methods.  Here I use class field theory to produce quartic or quadratic factors dividing $\hat H_{5,l}(x)$, instead of using direct computations in the function field of $E_5$, as was done in \cite{mor4} and \cite{mor5}.  \medskip

In Section 2 I prove a congruence modulo $l$ for the class equation $K_{5l}(X)=H_{-20l}(X)$ (for $l \equiv 1$ mod $4$) or the product of class equations $K_{5l}(X)=H_{-5l}(X) H_{-20l}(X)$ (for $l \equiv 3$ mod $4$), using the method of \cite{mor7}.  This requires $l$ to be a prime greater than $379$ (or one of a set $\mathcal{S}$ of $22$ specified primes less than $379$; see (\ref{eq:6})).  In Section 3 I show that each factor of the form $g(x)$ or $k(x)$ in Theorem 1.1 corresponds to an elliptic curve $E_5(b)$ for which $\mu=\sqrt{-5l}$ injects into the endomorphism ring $\textrm{End}(E_5(b))$, where $b$ is a root of $g(x)$ or $k(x)$ over $\mathbb{F}_l$.  The proof of this uses the formula in \cite{mor3} for points of order $5$ in $E_5(b)[5] - \langle (0,0) \rangle$. \medskip

This implies that all factors of $\hat H_{5,l}(x)$ of the form $g(x)$ or $k(x)$ in Theorem 1.1 arise from $j$-invariants which are roots of $K_{5l}(X)$ over $\mathbb{F}_l$.  In particular, the factors in Theorem 1.1 divide one of the polynomials
$$F_d(x) = x^{5\textsf{h}(-d)}(1-11x-x^2)^{\textsf{h}(-d)}H_{-d}\left(\frac{(x^4+12x^3+14x^2-12x+1)^3}{x^5(1-11x-x^2)}\right)$$
over $\mathbb{F}_l$, where $-d=-5l$ or $-d=-20l$ and $\textsf{h}(-d)$ is the class number of the order $\mathcal{O}=\textsf{R}_{-d}$ of discriminant $-d$ in $K=\mathbb{Q}(\sqrt{-5l})$.  This holds because of the formula for $\hat H_{5,l}(x)$ given above and the fact that
$$j(b) = j(E_5(b)) = \frac{(b^4+12b^3+14b^2-12b+1)^3}{b^5(1-11b-b^2)}$$
is the $j$-invariant of $E_5(b)$.  \medskip

In Section 4, I determine computationally how many factors in Theorem 1.1 arise from special factors $H_{-d}(X)$ of $K_{5l}(X)$ in Theorem 2.1, where $d$ equals $20$ or is one of the integers in the set
$$\mathfrak{T} = \{4, 11, 16, 19, 24, 36, 51, 64, 84, 91, 96, 99\}.$$
For $l>379$, the class equations corresponding to $d \in \mathfrak{T}$ are the factors which occur to the $4$-th power in $K_{5l}(X)$ mod $l$, when their roots are supersingular.  I show that each of the factors $H_{-d}$, for $d \in \mathfrak{T}$ and $d \neq 4$, yields $\textrm{deg}(H_{-d}(X))$ factors of the form $g(x)$, or twice that number of factors of the form $k(x)$, dividing $\hat H_{5,l}(x)$ over $\mathbb{F}_l$, depending on the congruence class of $l$ modulo $5$, for $l >379$ (or $l \in \mathcal{S}$).  \medskip

On the other hand, Theorem 2.1 shows that $K_{5l}(X)$ is also exactly divisible over $\mathbb{F}_l$ by a number of factors of the form $(X^2+a_iX+b_i)^2$.  I show in Sections 5 and 6 that each of these factors contributes exactly one factor of the form $g(x)$, or two factors of the form $k(x)$, to the factorization of $\hat H_{5,l}(x)$, depending on the congruence class of $l$ modulo $5$.  The fact that $(X^2+a_iX+b_i)^2$ exactly divides $K_{5l}(X)$ means that this factor can only divide one of the polynomials $H_{-5l}(X)$ or $H_{-20l}(X)$ modulo $l$.  The factors arising from $H_{-5l}(X)$ come about as follows.  Let $(5) =\wp_5^2$ and $(l) = \mathfrak{l}^2$ in $K = \mathbb{Q}(\sqrt{-5l})$.  The ray class field $\Sigma_{\wp_5}$ with conductor $\wp_5$ over $K$ is generated over $K$ by the fifth power of a value of the Rogers-Ramanujan continued fraction $r(\tau)$, where $\tau \in K$.  If $l \equiv 2, 3$ (mod $5$), the minimal polynomial of $r(\tau)^5$ over the decomposition field $L$ of $\mathfrak{l}$ in $\Sigma_{\wp_5}/K$ has the form $\tilde g(x) = x^4+a x^3+(11a+2)x^2-ax+1$, where $a = j_5^*(\tau)$ is a value of the modular function
$$j_5^*(\tau)= \left(\frac{\eta(\tau)}{\eta(5\tau)}\right)^6+22+125\left(\frac{\eta(5\tau)}{\eta(\tau)}\right)^6.$$
($\eta(\tau)$ is the Dedekind $\eta$-function.)  Reducing $\tilde g(x)$ modulo a prime divisor of $\mathfrak{l}$ in $L = K(j_5^*(\tau))$ yields an irreducible polynomial $g(x) \in \mathbb{F}_l[x]$ of the required form, which must divide $\hat H_{5,l}(x)$ because $l$ divides the discriminant of $K/\mathbb{Q}$, implying that the roots of $H_{-5l}(X)$ mod $l$ are supersingular $j$-invariants.  When $l \equiv 1, 4$ (mod $5$), $\tilde g(x)$ splits into two polynomials of the form $k(x)$ modulo a prime divisor of $\mathfrak{l}$ in the field $L = K(j_5^*(\tau))$ (a subfield of the decomposition field of $\mathfrak{l}$ in this case), yielding two of the required factors $k(x)$ modulo $l$.  When $l \equiv 3$ (mod $4$) there is a similar argument for the field $\Sigma_{\wp_5} \Omega_2$ ($\Omega_2$ is the ring class field of conductor $f=2$ over $K$), corresponding to factors arising from $H_{-20l}(X)$ in the above formula for $F_d(x)$.
\medskip

Using the fact that a root of the irreducible polynomial $H_{-d}(X)$ ($d = 5l$ or $20l$) generates a subfield over which the ring class field $\Omega_f$ ($f=1$ or $2$) has degree $2$, it is not hard to see that the factors $(X^2+a_iX+b_i)^2$ considered above correspond to distinct pairs of prime divisors of degree $2$ over $\mathfrak{l}$ in the extension $\Omega_f/K$.  Thus, the above argument yields roughly $\frac{1}{4} [\Omega_f:K] = \frac{1}{4}\textsf{h}(-d)$ (i.e. $\frac{1}{4}\textsf{h}(-d) + O(1)$) quartic factors of the form $g(x)$, when $l \equiv 2, 3$ (mod $5$).  When $l \equiv 1$ (mod $4$), $d=20l$ and $\textsf{h}(-20l) = h(-5l)$, so this gives the right number of quartic factors in Theorem 1.1A.  When $l \equiv 3$ (mod $8$), the total number of quartic factors is then roughly $\frac{1}{4}(\textsf{h}(-5l)+\textsf{h}(-20l)) = \frac{1}{2} h(-5l)$; and when $l \equiv 7$ (mod $8$), the total number of quartic factors is again roughly $\frac{1}{4}(\textsf{h}(-5l)+\textsf{h}(-20l)) = \frac{1}{4} (h(-5l)+3h(-5l)) = h(-5l)$.  These arguments show why the counts in Theorem 1.1A should be approximately true.  \medskip

To complete the argument we must show that each of the factors $(X^2+a_iX+b_i)^2$ of $K_{5l}(X)$ mod $l$ contributes exactly one factor of the form $g(x) = k_1(x) k_2(x)$, when $l \equiv 2, 3$ (mod $5$) and two factors $k(x)$ otherwise, to the factorization of $\hat H_{5,l}(x)$.  Writing the roots of $k_i(x) = x^2 + r_i x +s_i$ as $\alpha^5, \beta^5$, the condition $r_i = \varepsilon^5 (s_i-1)$ becomes
$$\alpha^5+\beta^5 = \varepsilon^5(1-\alpha^5 \beta^5).$$
It turns out that another factor of the form $g(x)$ would yield a different solution $\alpha' = M_1(\alpha), \beta' = M_2(\beta)$ of this equation, where $M_1, M_2$ are linear fractional maps contained in Fricke's normal form $G_{60} \cong A_5$ of the icosahedral group.  A sequence of suitable resultant calculations shows that this can only happen when the factor $g(x)$ arises not from $(X^2+a_iX+b_i)^2$ but from one of the factors $H_{-d}(X)$ with $d \in \mathfrak{T} - \{4,11,16,19\}$.  As shown before, these particular factors of $K_{5l}(X)$ do indeed yield multiple quartics $g(x)$, and the final argument shows these are the only factors of $K_{5l}(X)$ which do so.  All of this works for primes $l > 379$ or $l \in \mathcal{S}$.  The proof is then completed by the computations in Tables 6-9, which show that the formulas of Theorem 1.1 also hold for all primes $l$ for which $7 \le l \le 379$.  \medskip

As was shown in \cite[Theorem 7.1]{mor8}, the above theorem, together with the results of \cite[Theorems 1.1 and 1.3]{mor8}, implies the following result, which was conjectured by Nakaya \cite{na} (see his conjecture 5 for $N=5$). This theorem gives a formula for the number of linear factors over $\mathbb{F}_p$ of the supersingular polynomial $ss_p^{(5*)}(X)$ corresponding to the Fricke group $\Gamma_0^*(5)$, introduced by Koike and Sakai \cite{sa1}, \cite{sa2}.  As in \cite{na}, $L(p)$ denotes the number of supersingular $j$-invariants of elliptic curves which lie in the prime field $\mathbb{F}_p$.  This number was determined by Deuring \cite{d3} and is given by
$$L(p) = S(\mathbb{F}_p) = \begin{cases}
\frac{1}{2}h(-p), &\textrm{if} \ p \equiv 1 \ (\textrm{mod} \ 4),\\
2h(-p), &\textrm{if} \ p \equiv 3 \ (\textrm{mod} \ 8),\\
h(-p), &\textrm{if} \ p \equiv 7 \ (\textrm{mod} \ 8);
\end{cases}$$
where $h(-p)$ is the class number of the field $\mathbb{Q}(\sqrt{-p})$.  (See also \cite[p. 97]{brm}.)  \medskip

\noindent {\bf Theorem 1.2.}
{\it If $p>5$ is a prime, the number of linear factors of the supersingular polynomial $ss_p^{(5*)}(X)$ corresponding to the Fricke group $\Gamma_0^*(5)$ over $\mathbb{F}_p$ is given by the formula}
\begin{align*}
L^{(5*)}(p) & =  \ \frac{1}{2}\left(1+\left(\frac{-p}{5}\right)\right)L(p)\\
& +\ \frac{1}{8}\Big\{2+\left(1-\left(\frac{-1}{5p}\right)\right)\left(2+\left(\frac{-2}{5p}\right)\right)\Big\}h(-5p)\\
& = \begin{cases}
\frac{1}{4}\left(1+\left(\frac{p}{5}\right)\right)h(-p)+\frac{1}{4} h(-5p), &\textrm{if} \ p \equiv 1 \ \textrm{mod} \ 4;\\
\left(1+\left(\frac{p}{5}\right)\right)h(-p)+\frac{1}{2} h(-5p), &\textrm{if} \ p \equiv 3 \ \textrm{mod} \ 8;\\
\frac{1}{2}\left(1+\left(\frac{p}{5}\right)\right)h(-p)+ h(-5p), &\textrm{if} \ p \equiv 7 \ \textrm{mod} \ 8.
\end{cases}
\end{align*}

This result is analogous to formulas proved by Nakaya \cite{na} for the polynomials $ss_p^{(2*)}(X)$ and $ss_p^{(3*)}(X)$. \medskip

Nakaya \cite{na} has conjectured \cite[Conjectures 1, 6]{na} that the polynomial $ss_p^{(5*)}(X)$ has degree given by
$$\textrm{deg}(ss_p^{(5*)}(X)) = \frac{1}{4}\left(p-\left(\frac{-1}{p}\right)\right)+\frac{1}{2}\left(1-\left(\frac{-5}{p}\right)\right).$$
I give a proof of this conjecture in Section 7, using a parametrization discussed in \cite[Theorem 6.1]{mor8}.  A similar proof also establishes his conjecture for $\textrm{deg}(ss_p^{(7*)}(X))$.  \medskip

Nakaya has also conjectured that $ss_p^{(5*)}(X)$ is a product of linear factors (mod $p$) if and only if $p$ is one of the primes in the set $\{2, 3, 5, 7, 11,19\}$.  Curiously, these are exactly the prime divisors of the order of the Harada-Norton group $HN$ and of the Janko group $J_1$.  See \cite[Conjecture 2]{na}, \cite[Ch. 10]{con}.  Theorem 1.2 implies that for $p>5$, $ss_p^{(5*)}(X)$ is a product of linear factors (mod $p$) if and only if
$$L^{(5*)}(p) = \frac{1}{4}\left(p-\left(\frac{-1}{p}\right)\right)+\frac{1}{2}\left(1-\left(\frac{-5}{p}\right)\right).$$
Nakaya's Conjecture 2 (for $N=5$) can be obtained from this formula using a standard estimate of the class numbers $h(-p)$ and $h(-5p)$, along with a straightforward, if tedious, calculation.  (See Table 10 in Section 7.)  As he has discussed in \cite{na} in connection with a number of examples of the same phenomenon, this is an analogue of Ogg's observation that the $15$ primes for which the supersingular polynomial $ss_p(X)$ is a product of linear factors over $\mathbb{F}_p$ coincide with the prime factors of the order of the Monster group.  See the discussion in \cite{dgo} and \cite{do}.  See also \cite{bra3} and the related papers \cite{bra1} and \cite{bra2} of H. Brandt. \medskip

\section{Factorization of $K_{5p}(x)$ mod $p$}

Define $K_{5p}(X) = H_{-20p}(X)$ if $p \equiv 1$ (mod $4$) and $K_{5p}(X) = H_{-5p}(X) H_{-20p}(X)$ if $p \equiv 3$ (mod $4$).  The proof of Theorem 1.1 is based on the following congruence for $K_{5p}(X)$ modulo $p$.  \bigskip

\noindent {\bf Theorem 2.1.} {\it If $p > 379$ is a prime, then we have the factorization
$$K_{5p}(X) \equiv H_{-20}(X)^{2\epsilon_{20}} \prod_{d \in \mathfrak{T}}{H_{-d}(X)^{4\epsilon_d}} \prod_{i}{(X^2+a_iX+b_i)^2} \ (\textrm{mod} \ p),$$
where $H_{-d}(X)$ is the class polynomial for discriminant $-d$, $\mathfrak{T}$ is the set
$$\mathfrak{T} = \{4, 11, 16, 19, 24, 36, 51, 64, 84, 91, 96, 99\};$$
$\epsilon_d$ is defined by
\begin{align*}
\epsilon_{20} &=  \frac{1}{4}\left(1-\left(\frac{-20}{p}\right)\right)\left(1+\left(\frac{5}{p}\right)\right),\\
\epsilon_d &= \frac{1}{2}\left(1-\left(\frac{-d}{p}\right)\right), \ \textrm{if} \ d \in \{4,11,16,19\},\\
\epsilon_d &= \frac{1}{4}\left(1-\left(\frac{-d}{p}\right)\right)\left(1-\left(\frac{disc(H_{-d}(X))}{p}\right)\right),\ \textrm{if} \ d \in \{24,36,51,64,91,99\},
\end{align*}
\begin{align*}
\epsilon_{84} &=  \frac{1}{8}\left(1-\left(\frac{-84}{p}\right)\right)\left(1-\left(\frac{3}{p}\right)\right)\left(1-\left(\frac{7}{p}\right)\right),\\
\epsilon_{96} &=\frac{1}{8}\left(1-\left(\frac{-96}{p}\right)\right)\left(1-\left(\frac{2}{p}\right)\right)\left(1-\left(\frac{3}{p}\right)\right);
\end{align*}
and the polynomials $X^2+a_iX+b_i$ in the product $\prod_i{}$ are certain irreducible factors of the supersingular polynomial $ss_p(X)$ over $\mathbb{F}_p$ which are distinct from the factors in the product over $d \in \mathfrak{T}$.} \bigskip

For the sake of completeness, I note here the discriminants of the quadratic factors $H_{-d}(X)$, for $d \in \mathfrak{T}$:
\begin{align*}
\textrm{disc}(H_{-24}(X)) &= (2^{19}) 3^6 13^2 19^2,\\
\textrm{disc}(H_{-36}(X)) &= 2^{20} (3^3) 7^4 19^2 31^2,\\
\textrm{disc}(H_{-51}(X)) &= 2^{30} 3^6 7^4 (17) 31^2,\\
\textrm{disc}(H_{-64}(X)) &= (2^5) 3^{14} 7^4 11^4 19^2 59^2,\\
\textrm{disc}(H_{-91}(X)) &= 2^{32} 3^{12} 7^2 11^4 (13) 71^2,\\
\textrm{disc}(H_{-99}(X)) &= 2^{30} (3) 7^4 (11^3) 13^2 19^4 79^2.
\end{align*}
Factors in parentheses indicate nontrivial contributions to the Legendre symbols in Theorem 2.1.  Note that $d \in \mathfrak{T}$ if and only if $d \equiv \pm 1$ (mod $5$) and the prime ideal divisor $\wp_5$ of $5$ in $K = \mathbb{Q}(\sqrt{-d})$ satisfes $\wp_5^2 \sim1$ in the ring class group (mod $f$), where $-d = d_Kf^2$.  See \cite[Prop. 3.2]{mor2}. \medskip

\noindent {\it Proof of Theorem 2.1.} This is proved by the method of \cite{mor7} and an extended computer calculation. The proof will be divided into several parts.  \medskip

\noindent {\it Notation and method of proof.} \medskip

Let $\Phi_5(X,Y)=0$ denote the modular curve of level $5$.  This equation can be computed using the resultant
$$5^{15} \Phi_5(x,y) = \textrm{Res}_z((z^2+12z+16)^3+x(z+11),(z^2-228z+496)^3+y(z+11)^5).$$
We have by direct calculation on Maple that
\begin{align*}
\textrm{disc}_y(\Phi_5(x,y)) &= 5^5 x^4(x-1728)^4 \prod_{d \in \mathfrak{T}-\{4\}}{H_{-d}(x)^2},\\
\Phi_5(x,x) &= - (x^2-1264000x-681472000)(x-1728)^2(x+32^3)^2\\
 & \ \ \ \times (x-66^3)^2(x+96^3)^2\\
& = -H_{-20}(x) H_{-4}(x)^2 H_{-11}(x)^2 H_{-16}(x)^2 H_{-19}(x)^2.
\end{align*}
See also Fricke \cite[III, p.338]{fr} for this discriminant formula and Cox \cite[p. 263]{co} for the second formula.  Also, let $Q_5(u,v)$ denote the de-symmetrized form of $\Phi_5(X,Y)$, so that $Q_5(-x-y,xy) = \Phi_5(x,y)$.  This polynomial is:
\begin{align*}
Q_5(u,v)=& \ u^6 - 1963211489280u^5\\
& \ + (-246683410956v + 1284733132841424456253440)u^4\\
 & \ + (-2028551200v^2 - 128541798897012758937600v\\
 & \ - 280244777828439527804321565297868800)u^3\\
 & \ + (-4550940v^3 + 383083610766544859184v^2\\
 & \ - 192457939757860831020806056181760v\\
 & \ + 6692500042627997708487149415015068467200)u^2\\
 & \ + (-3720v^4 - 107878922099683200v^3\\
 & \ - 26898103232984020907026022400v^2\\
 & \ - 35714002250464310712293507636763033600v\\
 & \  - 53274330803424425450420160273356509151232000)u\\
 & \ - v^5 + 1666008466480v^4 - 441973132732967824498752v^3\\
 & \   + 5495857649359740948103830574202880v^2\\
 & \  - 277458457161876591676690089078008919883776v\\
 & \ + 2^{90} \cdot 3^{18} \cdot 5^3 \cdot11^9.
 \end{align*}
From \cite[Lemma 2.3]{mor7} we know that for $p > 20$ the irreducible factors of $K_{5p}(x)$ are the same as the supersingular factors which divide $\Phi_5(x^p,x)$, and their multiplicities are the same as their multiplicities in $\Phi_5(x^p,x)^2$:
$$K_{5p}(x) = \prod_{i}{q_i(x)^{e_i}}, \ \ \textrm{over} \ \ q_i(x) \mid \textrm{gcd}(ss_p(x),\Phi_5(x^p,x))$$ 
and $q_i(x)^{e_i} || \Phi_5(x^p,x)^2$.  For the proof of the theorem we shall calculate $F(x) = \Phi_5(x^p,x)$ and its first and second derivatives mod $p$ evaluated at a root $t$ of each supersingular factor $H_{-d}(X)$, for $d \in \mathfrak{T}$, to show that $H_{-d}(X)^4 \ || \ K_{5p}(X)$ mod $p$.  These are the only factors which can occur to a power greater than $2$ in $K_{5p}(X)$, because
$$\textrm{disc}_x(\Phi_5(x^p,x)) \mid \Delta = \textrm{disc}_y(\Phi_5(x,y)).$$
See \cite[Prop. 2.4]{mor7}. \medskip

\noindent {\it Linear factors.} \medskip

Roots of $K_{5p}(X)$ in $\mathbb{F}_p$ satisfy $x^p = x$, so the above remarks imply that linear factors can only come from the roots of $\Phi_5(x,x) \equiv 0$, i.e. factors of $H_{-d}(x)$ (mod $p$) for $d \in \{4, 11, 16, 19, 20\}$; and these factors, for $d \neq 20$, are supersingular in characteristic $p$ if and only if the corresponding $\epsilon_d =1$.  For $d=20$, the roots of
$$H_{-20}(x)=x^2-1264000x-681472000,$$
which lie in $\mathbb{Q}(\sqrt{5})$, are supersingular and contained in $\mathbb{F}_p$ if and only if $\epsilon_{20}=1$. We must check that the multiplicities of the linear factors are $2$, if $d=20$ and $4$, if $d \in \{4, 11, 16, 19\}$. \medskip

With $Q(u,v)=Q_5(u,v)$ let
$$Q_1 = \frac{\partial Q(u,v)}{\partial u}, \ \ Q_2 = \frac{\partial Q(u,v)}{\partial v},$$
denote the first partials of $Q$ and $Q_{ij}$ the second partials.  If
$$F(t) = \Phi_5(t^p,t)=Q(-t^p-t,t^{p+1}),$$
then in characteristic $p$,
\begin{align*}
F'(t) & = -Q_1(-t^p-t, t^{p+1}) + t^p Q_2(-t^p-t, t^{p+1}),\\
F''(t) & = Q_{11}(-t^p-t, t^{p+1}) - 2t^pQ_{12}(-t^p-t, t^{p+1}) + t^{2p}Q_{22}(-t^p-t, t^{p+1}).
\end{align*}
For the roots $t$ of the linear factors we have
\begin{align}
\label{eq:1}
F(t) & =Q(-2t,t^2), \ \ F'(t) =-Q_1(-2t,t^2)+tQ_2(-2t,t^2),\\
\notag F''(t) & = Q_{11}(-2t,t^2)-2tQ_{12}(-2t,t^2) + t^2Q_{22}(-2t,t^2).
\end{align}
The following values can be checked on Maple:
\begin{align*}
& F(1728)=F'(1728)=0, \ F''(1728)=2^{41} 3^{24} 5^2 7^8 11^4 19^4;\\
& F(-32^3)=F'(-32^3)=0, \ F''(-32^3)=-2^{63} \cdot 5\cdot 7^8 11^2 13^4 17^2 19^3 43^2;\\
& F(66^3)=F'(66^3)=0, \ F''(66^3)=2^{21} 3^{26} \cdot 5 \cdot 7^8 11^3 19^3 \cdot 31 \cdot 43^2 67^2 \cdot 71 \cdot 79;\\
& F(-96^3)=F'(-96^3)=0, \ F''(-96^3)=-2^{63} 3^{26} \cdot 5 \cdot 13^4 19^2 \cdot 31 \cdot 59 \cdot 67^2 \cdot 79.
\end{align*}
These values show that the multiplicity of $H_{-d}(x)$ in $\Phi_5(x^p,x)$ over $\mathbb{F}_p$ is $2$, for $d \in \{4, 11, 16, 19\}$ and $p > 79$, and therefore the multiplicities of the $q_i(x)=H_{-d}(x)$ in $K_{5p}(x)$ for these $d$ are $e_i=4$, when they occur (for $p>79$).  Moreover, the linear factors corresponding to these values of $d$ are distinct (mod $p$), for $p > 67$. \medskip

For $d=20$, we use (\ref{eq:1}) to evaluate $F(t)$ and its derivative at $t = 632000 + 282880 \sqrt{5}$, which is a root of $H_{-20}(x)$.  We find that $F(t)=0$ in characteristic $0$, but that $F'(t) = A+B \sqrt{5}$, with
\begin{align*}
A = &-2^{56} \cdot 3 \cdot 5^5 \cdot 7 \cdot 11^7 \cdot 13^5 \cdot 17^3 \cdot 19^4 \cdot 31^3 \cdot 79^2\cdot 919\\
B = & -2^{54} \cdot 5 \cdot 11^6 \cdot 13^5 \cdot 17^3 \cdot 19^4 \cdot 29 \cdot 31^2 \cdot 79^2 \cdot 467 \cdot 543287;
\end{align*}
giving that
$$A^2-5B^2 = -2^{108} \cdot 5^3 \cdot 11^{12} \cdot 13^{10} \cdot 17^6 \cdot 19^{10} \cdot 31^4 \cdot 59^2 \cdot 71^2 \cdot 79^4.$$
This calculation shows that the linear factors $(x-t), (x-t')$ of $H_{-20}(x)$ divide $K_{5p}(x)$ to the second power, when they occur and lie in $\mathbb{F}_p[x]$ (for $p>79$).  Factoring $H_{-20}(t)$ for $t = 12^3, -32^3, 66^3, -96^3$ shows that none of the linear factors mentioned above coincides with a factor of $H_{-20}(x)$ (mod $p$), when $p > 79$.  This completes the argument for the linear factors of $K_{5p}(x)$.  \medskip

\noindent {\it Quadratic factors.} \medskip

Next, the class equations $H_{-d}(x)$ are quadratic for the six values of $d \in \{24,36,51,64,91,99\}$.  Letting $t$ be a root of $H_{-d}(x)=x^2+u x+v$, for one of these $d$, we verify that $F(t)=Q(u,v)=0$ and $F'(t)=Q_1(u,v)=Q_2(u,v)=0$ in characteristic $0$, but that $F''(t) \not \equiv 0$ in $\mathbb{F}_p$, for $p> 379$.  To check the latter note that $F''(t) = D_1(u,v)-t^p D_2(u,v)$, where
$$D_1(u,v)= Q_{11}(u,v)-vQ_{22}(u,v), \ \ D_2(u,v) = 2Q_{12}(u,v)+uQ_{22}(u,v).$$
(See \cite[eq. (2.8)]{mor7}.)  Using that $1, t^p$ are independent over $\mathbb{F}_p$, it follows that $F''(t) \equiv 0$ mod $p$ if and only if $p \mid D_1(u,v)$ and $p \mid D_2(u,v)$.  We check that $F''(t) \not \equiv 0$ (mod $p$) for the six values of $d$ above by evaluating $D_1$ and $D_2$ at the coefficients of the polynomials $H_{-d}(x)$.  Recall now that
\begin{align*}
H_{-24}(x)&=x^2-4834944x+14670139392,\\
H_{-36}(x)&=x^2-153542016x-1790957481984,\\
H_{-51}(x)&=x^2+5541101568x+6262062317568,\\
H_{-64}(x)&=x^2-82226316240x-7367066619912,\\
H_{-91}(x)&=x^2+10359073013760x-3845689020776448,\\
H_{-99}(x)&=x^2+37616060956672x-56171326053810176.
\end{align*}
The values of $\textrm{gcd}(D_1,D_2)$ for these values of $d$ are contained in Table 1.  The factorizations in the table show that the quadratic polynomials $H_{-d}(x)$ divide $K_{5p}(x)$ to exactly the fourth power (mod $p$), when they occur, for $p > 379$.  Furthermore, factoring the coefficients of the differences $H_{-d_1}(x)-H_{-d_2}(x)$ shows that these class polynomials are all distinct (mod $p$), when $p > 307$. \medskip

\begin{table}
  \centering 
  \caption{Calculating $F''(t) \neq 0$.}\label{tab:0}

\noindent \begin{tabular}{|c|l|c|}
\hline
   &   \\
$d$	&   \ \ \ \ $\textrm{gcd}(D_1,D_2)$   \\
\hline
  &    \\
 24  &  $2^{36} \cdot 3^{18} \cdot 5 \cdot 13^3 \cdot 19^2 \cdot 37 \cdot 43^2 \cdot 61 \cdot 67 \cdot 109$  \\
  &  \\
 36  &  $2^{39} \cdot 3^8 \cdot 5 \cdot 7^6 \cdot 19^2 \cdot 43^2 \cdot 67 \cdot 79 \cdot 127 \cdot 139 \cdot 151 \cdot 163$  \\
  &  \\
 51 & $2^{46} \cdot 3^{19} \cdot 5 \cdot 7^7 \cdot 17 \cdot 37 \cdot 61 \cdot 79 \cdot 139 \cdot 163 \cdot 211$ \\
  &  \\
 64 & $2^{11} \cdot 3^{20} \cdot 5 \cdot 7^6 \cdot11^3 \cdot 19^2 \cdot 43^2 \cdot 67 \cdot 139 \cdot 163 \cdot 211 \cdot 283 \cdot 307$ \\
  &  \\
 91 & $2^{50} \cdot 3^{18} \cdot 5 \cdot 7^4 \cdot 11^2 \cdot 13^2 \cdot 37 \cdot 61 \cdot 67 \cdot 109 \cdot 151 \cdot 163 \cdot 331 \cdot 379$ \\
  & \\
 99 & $2^{46} \cdot 5 \cdot 7^6 \cdot 11 \cdot 13^3 \cdot 19^2 \cdot 29 \cdot 41 \cdot 43^2 \cdot 61 \cdot 109 \cdot 127 \cdot 139 \cdot 211 \cdot 283 \cdot 307$ \\
  \hline
\end{tabular}
\end{table}

\noindent {\it Quadratic factors}, $d=84$. \medskip

The only other quadratic factors that can occur in the factorization of $K_{5p}(x)$ to a power higher than the second are quadratic factors of the quartic polynomials $H_{-84}(x)$ and $H_{-96}(x)$, by the formula for the discriminant of $\Phi_5(x,y)$.  We start with the polynomial
\begin{align*}
H_{-84}(x) = & \ x^4 - 3196800946944x^3 - 5663679223085309952x^2\\
& + 88821246589810089394176x - 5133201653210986057826304,\\
= & \ (x^2 - 1598400473472x + 92704725504000)^2\\
& -2^{18}3^9 7^3 13^2 29^2 (3187x - 184896)^2,
\end{align*}
whose discriminant is
$$\textrm{disc}(H_{-84}(x))=2^{116}3^{44}7^{14}13^{12}29^4 43^2 53^2 61^2 67^2 73^2 79^2.$$
With some calculation using the above expression, it is straightforward to check that the splitting field of $H_{-84}(x)$ is $\mathbb{Q}(\sqrt{3},\sqrt{7})$.  If $\left(\frac{3}{p}\right)=+1$ and $\left(\frac{7}{p}\right)=-1$, where $p > 79$, then $H_{-84}(x)$ factors into two irreducible quadratics over $\mathbb{F}_p$.  One of these factors is
\begin{align*}
q_1(x)=x^2+ux+v = & \ x^2 + (922836934656\sqrt{3} - 1598400473472)x\\
& + 1649310419952599040\sqrt{3} - 2856689444809764864,
 \end{align*}
but
\begin{align}
\notag \textrm{Norm}_\mathbb{Q}(Q(u,v)) = & \ 2^{108}3^{62}7^{12}13^{12} 29^6 43^2 47^2 53^4 61^3 73^2 97^3\\
& \times 157 \cdot 181 \cdot 229 \cdot 241 \cdot 313 \cdot 349 \cdot 397 \cdot 409.
\label{eq:2}
\end{align}
Hence, $q_1(x)$ and its conjugate $\tilde q_1(x)$ over $\mathbb{Q}(\sqrt{3})$ do not divide $\Phi_5(x^p,x)$ (mod $p$), for $p > 409$.  Furthermore,
$$\textrm{gcd}(N(Q_1),N(Q_2))=2^{76}3^{43}7^8 13^8 29^4 \cdot 47 \cdot 53^2 61^2 \cdot 73 \cdot 97,$$
where $N(a)$ denotes the norm to $\mathbb{Q}$, so $q_1(x)$ can occur only to the second power in $K_{5p}(x)$ and can be absorbed into the final product of the theorem, when it or $\tilde q_1(x)$ occurs, for $97<p \le 409$.  In Section 6 we will call these factors, as well as similar factors occurring below, {\it sporadic} factors.  (See the analogous arguments below using equations (\ref{eq:3}), (\ref{eq:4}), (\ref{eq:5}).)
\medskip

Similarly, if $\left(\frac{3}{p}\right)=-1$ and $\left(\frac{7}{p}\right)=+1$, $H_{-84}(x)$ has the factor
\begin{align*}
q_2(x) =x^2+ux+v = & \ x^2 + (604139268096\sqrt{7} - 1598400473472)x\\
 &- 9357315081633792 \sqrt{7} + 24757128541605888.
 \end{align*}
 In this case
 \begin{align}
\notag \textrm{Norm}_\mathbb{Q}(Q(u,v)) = & \ 2^{108}3^{36}7^{15}13^{12}29^4 43^2 47^2 \cdot 53 \cdot 59^2 61^4 73^4 83^2\\
\times \ 113 \ \cdot &\ 131 \cdot 137 \cdot 149 \cdot 197 \cdot 233^2 \cdot 281 \cdot 317 \cdot 389 \cdot 401.
\label{eq:3}
\end{align}
Hence, $q_2(x)$ and $\tilde q_2(x)$ do not divide $\Phi(x^5,x)$ over $\mathbb{F}_p$, for $p>401$.  We have that
$$\textrm{gcd}(N(Q_1),N(Q_2))=2^{76}3^{30}7^{10} 13^8 \cdot 47 \cdot 61^2 73^2 \cdot 83,$$
and so $q_2(x)$ and $\tilde q_2(x)$ occur only to the second power in $K_{5p}(x)$, when either occurs, for $83 < p \le 401$.
\medskip

On the other hand, if $\left(\frac{3}{p}\right)=-1$ and $\left(\frac{7}{p}\right)=-1$, $H_{-84}(x)$ has the factor
\begin{align*}
q_3(x) =x^2+ux+v = & \ x^2 + (348799965696\sqrt{21} - 1598400473472)x\\
 &-20235870240768\sqrt{21} + 92704725504000,
 \end{align*}
 and in this case $Q(u,v)=Q_1(u,v)=Q_2(u,v)=0$ in characteristic $0$.  Since
\begin{align*}
\textrm{gcd}(N(D_1),N(D_2))= & \ 2^{70}3^{30}5^2 7^8 13^6 29^2 43^3 67^2\\
& \times 79 \cdot 127 \cdot 151 \cdot 163 \cdot 211^2 \cdot 331 \cdot 379,
\end{align*}
it follows that $q_3(x)$ and its conjugate $\tilde q_3(x)$ divide $K_{5p}(x)$ over $\mathbb{F}_p$, with multiplicity $4$, exactly when $\epsilon_{84} = 1$, for primes $p > 379$.  Moreover, the largest prime factor of any resultant $Res(H_{-84}(x),H_{-d}(x))$, for the values $d \in \{24, 36, 51, 64, 91, 99\}$, for which $H_{-84}(x)$ (mod $p$) is a product of two irreducible quadratics, is $379$, so the factors $q_3$ and $\tilde q_3$ are distinct from the quadratic factors found above, for $p >379$.  This establishes the contribution of $H_{-84}(x)$ to the factorization of $K_{5p}(x)$.  \medskip

\noindent {\it Quadratic factors}, $d=96$. \medskip

A similar analysis applies to the polynomial
\begin{align*}
H_{-96}(x) = & \ x^4 - 23340144296736x^3 + 670421055192156288x^2\\
& + 447805364111967209472x - 984163224549635621646336\\
= & \ (x^2 - 11670072148368x + 10900447400376000)^2\\
& - 2^93^{13}13^217^241^261^2(739x - 690264)^2,
\end{align*}
whose discriminant is
$$\textrm{disc}(H_{-96}(x))=2^{56}3^{46}13^{12}17^{12}19^6 23^2 37^2 41^4 43^2 61^4 67^2 89^2.$$
For $\left(\frac{2}{p}\right)=+1$ and $\left(\frac{3}{p}\right)=-1$, the polynomial
\begin{align*}
x^2+u x+v =&x^2 + (8251987131648\sqrt{2} - 11670072148368)x\\
& + 12701433452887296\sqrt{2} - 17962539423257664
\end{align*}
is a factor of $H_{-96}(x)$ (mod $p$) and
\begin{align}
\notag \textrm{Norm}_\mathbb{Q}(Q(u,v))=&2^{54}3^{36}13^{12}17^819^4 23^{10}37^4 41^2 43^2 47^2 61^6 67^2\\
& \times  89 \cdot 113^2 137^2 139^2 257 \cdot 281 \cdot 353 \cdot 401 \cdot 449.
\label{eq:4}
\end{align}
Moreover, 
$$\textrm{gcd}(N(Q_1),N(Q_2))=2^{40}3^{30}13^8 17^4 23^5 37^2 \cdot 47 \cdot 61^4 \cdot137.$$
Thus, the above factor only occurs to the second power in $K_{5p}(x)$, when it or its conjugate occurs, for $137<p \le 449$, so it can be absorbed into the final product in the congruence of the theorem. \medskip

For $\left(\frac{2}{p}\right)=-1$ and $\left(\frac{3}{p}\right)=+1$, the polynomial
\begin{align*}
x^2+u x+v =&x^2 + (6737719296672\sqrt{3} - 11670072148368)x \\
&- 197611189074074880\sqrt{3} + 342272619618959808
\end{align*}
is a factor of $H_{-96}(x)$ (mod $p$) and
\begin{align}
\notag \textrm{Norm}_\mathbb{Q}(Q(u,v))=&-2^{54}3^{63}13^{12}17^8 19^4 23^2 37^5 41^6 43^2 61^2 67^2 89^4\\
& \times 109^3\cdot 229 \cdot 277^2 \cdot 349 \cdot 373 \cdot 397 \cdot 421.
\label{eq:5}
\end{align}
In this case,
$$\textrm{gcd}(N(Q_1),N(Q_2))=2^{40}3^{49}13^817^437^241^4 \cdot 71 \cdot 89^2 \cdot 109,$$
so these factors occur only to the second power in $K_{5p}(x)$ when they occur, for $109<p \le 421$.
Finally, the polynomial
\begin{align*}
x^2 + ux +v = & x^2 - (4764286992816\sqrt{6} + 11670072148368)x\\
& + 4450089034924416\sqrt{6} + 10900447400376000
\end{align*}
is a factor of $H_{-96}(x)$ when $\left(\frac{2}{p}\right)=\left(\frac{3}{p}\right)=-1$, and $Q(u,v)=Q_1(u,v)=Q_2(u,v)=0$, while
\begin{align*}
\textrm{gcd}(N(D_1),N(D_2))= & \ 2^{41}3^{37}5^2 13^6 17^2 19^4 41^2 43^3 61^2 67^2 139^2 \\
& \times 163 \cdot 211 \cdot 283 \cdot 307 \cdot 331 \cdot 379.
\end{align*}
Thus, the factors of $H_{-96}(x)$ divide $K_{5p}(x)$ to the $4$-th power when $\epsilon_{96}=1$ and $p>379$.  The largest prime dividing a resultant $Res(H_{-96}(x),H_{-d}(x))$, for $d \in \{24, 36, 51, 64, 84, 91, 99\}$, modulo which $H_{-96}(x)$ factors as a product of two irreducible quadratics, is $379$.  This completes the discussion of the factors in the first product.
\medskip

The remaining irreducible quadratic factors $q_i(x)=x^2+a_ix+b_i$ are exactly the factors of $ss_p(x)$, distinct from the factors of $H_{-d}(x)$ ($d \in \mathfrak{T}$), for which $Q_5(a_i,b_i)=Q(a_i,b_i) \equiv 0$ (mod $p$), by \cite[Theorem 3.1]{mor7}, and their multiplicities in $K_{5p}(x)$ are exactly $2$.  This completes the proof of Theorem 2.1.  $\square$
\bigskip

For $p=379$, the polynomials $H_{-91}(x), H_{-84}(x), H_{-96}(x)$ all have the factor $H_{-91}(x) \equiv x^2+114x+51$ in common (mod $379$), and this factor occurs to the power $6$ in the factorization of $K_{5p}(x)$.  Thus, the factorization formula in Theorem 2.1 does not hold for $p=379$ (barely!).  By factoring the supersingular polynomial $ss_{379}(x)$ (mod $379$) and comparing with the factors and multiplicities of $\Phi_5(x^{379},x)$ (mod $379$), it can be checked that
\begin{align*}
& K_{5 \cdot 379}(x) = H_{-5 \cdot 379}(x) H_{-4 \cdot 5 \cdot 379}(x) \equiv \\
& \ \ \ \ (x + 163)^2(x + 181)^2(x+150)^4(x+165)^4(x+167)^4 \\
& \times (x^2 + 338x + 303)^4(x^2 + 359x + 73)^4(x^2 + 288x + 354)^4 \\
& \times (x^2 + 180x + 346)^4(x^2 + 114x + 51)^6 (x^2 + 47x + 352)^4 (x^2 + 23x + 346)^4 \\
& \times (x^2 + 68x + 125)^2(x^2 + 191x + 240)^2(x^2 + 320x + 244)^2(x^2 + 152x + 232)^2\\
& \times (x^2 + 57x + 374)^2 \ \ \ (\textrm{mod} \ 379).
\end{align*}
This agrees with the fact that $h(-5 \cdot 379)=h(-20 \cdot 379) = 48$.  Therefore, the condition $p >379$ in Theorem 2.1 is sharp. \medskip

Moreover, there are $22$ primes $p<379$ which do not occur in any of the factorizations in the proof of Theorem 2.1, and which do not divide the differences $H_{-d_1}(t)-H_{-d_2}(t)$ for $d_1, d_2 \in \{4,11,16,19,20\}$ or for $d_1,d_2 \in \{24, 36, 51, 64, 91,99\}$; or the remainders on dividing $H_{-84}(t)$ or $H_{-96}(t)$ by $H_{-d}(t)$, for $d \in \{24, 36, 51, 64, 91,99\}$; or the resultant $\textrm{Res}(H_{-84}(t),H_{-96}(t))$, for which $\left(\frac{3}{p}\right) = -1$.  These are the primes in the set
\begin{align}
\notag \mathcal{S}=\{&101, 103, 107, 167, 173, 179, 191, 193, 199, 223, 227,\\
& 239, 251, 263, 269, 271, 293, 311, 337, 347, 359, 367\}.
\label{eq:6}
\end{align}
For these primes all the arguments in the proof are valid.  This implies the following. \bigskip

\noindent {\bf Corollary 2.2.} {\it The assertion of Theorem 2.1 also holds for all $22$ primes in the set $\mathcal{S}$.}
\medskip

Equating degrees in the congruence of Theorem 2.1 yields the following formula.  Let
$$a_p = 1+\frac{1}{2}\left(1-\left(\frac{-1}{p}\right)\right) \left(2+\left(\frac{2}{p}\right)\right);$$
so that $a_p=1, 2, 4$ according as $p \equiv 1$ mod $4$, or $p \equiv 3, 7$ mod $8$. \bigskip

\noindent {\bf Theorem 2.3.} {\it For primes $p \in \mathcal{S}$ and for $p > 379$, the following formula holds:
$$a_p h(-5p) = 4\epsilon_{20}+\sum_{d \in \mathfrak{T}}{4\epsilon_d \textrm{deg}(H_{-d}(X))}+ 4N_p,$$
where $N_p$ is the number of irreducible quadratic factors $X^2+a_iX+b_i$ of $J_p(X)$, not dividing any of the 
factors $H_{-d}(X)^{\epsilon_d}$ (mod $p$) in Theorem 2.1, for which $Q_5(a_i,b_i) = 0$ (mod $p$).}
\bigskip

The factors $X^2+a_iX+b_i$ of $J_p(X)$ in this theorem, whose count is $N_p$, occur to only the first power in $\Phi_5(X^p,X)$ over $\mathbb{F}_p$.  See \cite[pp. 78, 83, 91, 92]{mor7}.  This yields an easy way of distinguishing them from the other factors of $J_p(X)$ for which $Q_5(a_i,b_i)=0$ in $\mathbb{F}_p$.  \medskip

\section{The endomorphism $\mu=\sqrt{-5l}$}

The next step in the proof is to show that for roots $b$ of the factors $g(x)$ and $k(x)$ in Theorem 1.1, the supersingular elliptic curve $E_5(b)$ has an endomorphism $\mu$ with $\mu^2 = -5l$, i.e., that $\sqrt{-5l}$ injects into the quaternion algebra $\textrm{End}(E_5(b))$.  To prepare for this, we solve the equation $g(x)=x^4+ax^3+(11a+2)x^2-ax+1=0$ algebraically, using the cubic resolvent for $g(x-\frac{a}{4})$.  (See \cite[pp. 194-196]{vdw}.) \medskip

\noindent {\bf Lemma 3.1.} {\it If $g(x) = x^4+ax^3+(11a+2)x^2-ax+1$, the roots of the cubic resolvent of $g(x-\frac{a}{4})$ are
\begin{align*}
\Theta_1 & = \frac{-1}{4}(a^2-44a-16),\\
\Theta_2 & = -a\left(\frac{a}{4}-\frac{11+5\sqrt{5}}{2}\right)=-a\left(\frac{a}{4}+\bar \varepsilon^5\right),\\
\Theta_3 & = -a\left(\frac{a}{4}-\frac{11-5\sqrt{5}}{2}\right)=-a\left(\frac{a}{4}+ \varepsilon^5\right).
\end{align*}
Moreover, $\Theta_2\Theta_3=\frac{a^2}{16}(a^2-44a-16)=-\frac{a^2}{4}\Theta_1$.  The roots of $g(x)=0$ are:} 
\begin{align*}
\rho_1 & :=\frac{-a}{4}- \frac{\sqrt{-\Theta_1}}{2}+ \frac{\sqrt{-\Theta_2}}{2} +  \frac{\sqrt{-\Theta_3}}{2},\\
\rho_2 & :=\frac{-a}{4}+ \frac{\sqrt{-\Theta_1}}{2}+ \frac{\sqrt{-\Theta_2}}{2} -  \frac{\sqrt{-\Theta_3}}{2},\\
\rho_3 & :=\frac{-a}{4}- \frac{\sqrt{-\Theta_1}}{2}- \frac{\sqrt{-\Theta_2}}{2} -  \frac{\sqrt{-\Theta_3}}{2},\\
\rho_4 & :=\frac{-a}{4}+ \frac{\sqrt{-\Theta_1}}{2}- \frac{\sqrt{-\Theta_2}}{2} +  \frac{\sqrt{-\Theta_3}}{2}.
\end{align*}

Straightforward calculation shows that
$$-(\rho_1+\rho_4)=\varepsilon^5(\rho_1 \rho_4-1), \ \ -(\rho_2+\rho_3)=\varepsilon^5(\rho_2 \rho_3-1).$$
It follows that $(x-\rho_1)(x-\rho_4)$ and $(x-\rho_2)(x-\rho_3)$ are polynomials of the form $x^2+rx+s$, with $r=\varepsilon^5(s-1)$, as in Theorem 1.1B.  Further,
$$-(\rho_1+\rho_2)=\bar \varepsilon^5(\rho_1 \rho_2-1), \ \ -(\rho_3+\rho_4)=\bar \varepsilon^5(\rho_3 \rho_4-1),$$
so that $(x-\rho_1)(x-\rho_2)$ and $(x-\rho_3)(x-\rho_4)$ have the form $x^2+rx+s$, with $r=\bar \varepsilon^5(s-1)$. \medskip

We also note that $\rho_3=-1/\rho_1$ and $\rho_4=-1/\rho_2$.  Thus, the factors of the form $k(x)$ in Theorem 1.1B and C come in pairs of factors, whose product has the form $g(x)$.  Note also that the roots of $g(x)$ are invariant under $\tau(b)=\frac{-b+\varepsilon^5}{\varepsilon^5b+1}$, since
\begin{equation}
(\varepsilon^5 x+1)^4 g(\tau(x)) = 5^3 \varepsilon^{10} g(x), \ \ \tau(x)=\frac{-x+\varepsilon^5}{\varepsilon^5 x+1}.
\label{eq:7}
\end{equation}
Also, $-1/\tau(b)=\bar \tau(b) = \frac{-b+\bar \varepsilon^5}{\bar \varepsilon^5b+1}$. \medskip

Now assume $b$ is a root of the factor $g(x)=x^4+ax^3+(11a+2)x^2-ax+1$ of $\hat H_{5,l}(x)=0$ over $\mathbb{F}_l$.  Then $E_5=E_5(b)$ is supersingular in characteristic $l$.  The calculations of \cite{mor} imply that $E_5(b)$ is isogenous to $E_5(\tau(b))$ by an isogeny $\phi_b$ of degree $5$, for the following reason.  By \cite[p. 259]{mor}, there is an isogeny $\psi:E_5 \rightarrow E_{5,5}$ defined over $\mathbb{F}_l(b)$, where $E_{5,5}=E_{5,5}(b)$ is the curve
\begin{align*}
E_{5,5}(b): \ Y^2+(1+b)XY+5bY=X^3 &+7bX^2+6(b^3+b^2-b)X\\
& + b^5+b^4-10b^3-29b^2-b,
\end{align*}
for which
$$X(\psi(P)) = \frac{b^4+(3b^3+b^4)x+(3b^2+b^3)x^2+(b-b^2-b^3)x^3+x^5}{x^2(x+b)^2}, \ \ x=X(P).$$
Furthermore,  $E_{5,5}(b) \cong E_5(\tau(b))$ by an isomorphism $\iota$, since these two curves have the same $j$-invariant, namely
$$j(E_{5,5}(b))=\frac{(b^4-228b^3+494b^2+228b+1)^3}{b(1-11b-b^2)^5}.$$
The $X$-coordinate of $\iota(Q)$ (for $Q \in E_{5,5}(b)$) is defined over $\mathbb{F}_{l^2}(b)$, since it is given by $\iota(X_1,Y_1)=(X_2,Y_2)$, where
$$\iota(X_1)=X_2= \lambda^2 X_1+\lambda^2 \frac{b^2+30b+1}{12}-\frac{\tau(b)^2+6\tau(b)+1}{12},$$
with $\displaystyle \lambda^2 = \frac{\sqrt{5}\bar \varepsilon^5}{(b-\bar \varepsilon^5)^2}$.  Composing the isogeny $\psi$ with this isomorphism $\iota$ gives the isogeny $\phi_b =\iota \circ \psi$.  \bigskip

\noindent {\bf Lemma 3.2.}  {\it Let $\phi_b: E_5(b) \rightarrow E_{5}(\tau(b))$ be the isogeny defined above, and let $\phi_{\tau(b)}: E_5(\tau(b)) \rightarrow E_{5}(b)$ be the isogeny obtained by replacing $b$ in the formulas for $\phi_b$ by $\tau(b)$ (leaving $\sqrt{5}$ fixed).  Then $\phi_{\tau(b)} \circ \phi_b = \alpha \circ [5]$, where $\alpha \in \textrm{Aut}(E_5(b))$ and $[5]$ is the multiplication-by-$5$ map on $E_5(b)$.  Thus, if $j(E_5(b)) \neq 0, 1728$, we have $\alpha = \pm 1$ and $\phi_{\tau(b)} = \pm \hat \phi_b$, where $\hat \phi_b$ is the dual isogeny of $\phi_b$.} \medskip

\noindent {\it Proof.}  Since $\tau(\tau(b)) = b$, it is clear that $\phi_{\tau(b)} \circ \phi_b$ is an isogeny from $E_5(b)$ to itself whose kernel contains the group $\langle (0,0) \rangle$.  To prove the lemma we must show that the kernel is actually $E_5(b)[5]$ and not a cyclic subgroup of $E_5(b)$ of order $25$.  For this we use the following formula from \cite{mor3} for the $X$-coordinate of a point $P \in E_5(b)$ of order $5$, which does not lie in $\langle (0,0) \rangle$:
\begin{align*}
X(P)=\frac{-\varepsilon^4}{2}\frac{(- 2u^2+(1+\sqrt{5})u - 3\sqrt{5} - 7)(2u^2 + (2\sqrt{5} + 4)u + 3\sqrt{5} + 7)}{(- 2u^2+(\sqrt{5}+1)u - 2)(u + 1)^2},
\end{align*}
where
$$u^5 = -\frac{b-\bar \varepsilon^5}{b-\varepsilon^5}.$$
A calculation on Maple shows that
$$X(\psi(P))= \frac{-5+\sqrt{5}}{10}(b^2+\varepsilon^4 b+\bar \varepsilon^2).$$
Then the above formula for $\iota$ gives $X_2 = \iota(X(\psi(P))) = 0$.  Hence, the point $P$ maps to $\phi_b(P)=\pm (0,0)$ on $E_5(\tau(b))$.  Now the kernel of $\psi$ is the group $\langle (0,0) \rangle$ on $E_5(b)$, so $\textrm{ker}(\phi_b) = \langle (0,0) \rangle$, whence it follows that $\textrm{ker}(\phi_{\tau(b)}) = \langle (0,0) \rangle$ on $E_5(\tau(b))$.  But the point $P \notin \langle (0,0) \rangle$, so $(0,0)$ and $P$ generate $E_5(b)[5]$.  Since $\phi_{\tau(b)} \circ \phi_b(P) = \phi_{\tau(b)}(\pm(0,0)) = O$ on $E_5(b)$, we have $\textrm{ker}(\phi_{\tau(b)} \circ \phi) = E_5(b)[5]$.  Since this is also the kernel of the multiplication-by-$5$ map on $E_5(b)$, it follows that $\phi_{\tau(b)} \circ \phi_b = \alpha \circ [5]$ in $\textrm{End}(E_5(b))$, for some automorphism $\alpha$ of $E_5(b)$.  Since $l >5$, the only possibility for $j(E_5) \neq 0, 1728$ is $\alpha =\pm [1]$.  This proves the lemma. (See \cite[pp. 73-74, 103]{sil2}.)   $\square$  \bigskip

\noindent {\bf Remark.} The minimal polynomials of the only values of $b$, for which $j(E_5(b)) = 0, 1728$, divide the respective polynomials
\begin{align*}
c_4(x) = & \ x^4+12x^3+14x^2-12x+1,\\
c_6(x) =  & \ -(x^2 + 1)(x^4 + 18x^3 + 74x^2 - 18x + 1).
\end{align*}
The quartic factor of $c_6(x)$ only has the form of $g(x)$ in characteristic $l > 5$ if $l=7$, and $c_4(x)$ never has the form of $g(x)$, for characteristic $l > 5$.  Note also that the discriminants of these quartic polynomials are only divisible by the primes $2, 3, 5$, so neither can be the square of a factor of the form $k(x)$, for $p > 5$.  \medskip

\noindent {\bf Lemma 3.3.} {\it Assume $E_5(b)$ is supersingular in characteristic $l$.} \smallskip

\noindent (a) {\it If $l \equiv 1, 4$ mod $5$, the multiplication-by-$l$ map $[l] \in End(E_5(b))$ acts on points $P \in E_5(b)$ by}
$$[l]P = \pm(X(P)^{l^2}, Y(P)^{l^2}) = \pm P^{l^2}.$$
\noindent (b) {\it If $l \equiv 2, 3$ mod $5$, we have}
$$[l^2]P = -(X(P)^{l^4}, Y(P)^{l^4}) = -P^{l^4}.$$

\noindent {\it Proof.}
\noindent (a) This is proved in the course of proving Theorem 6.1(i) in \cite[pp. 263-265]{mor}. (Take $m = 1$ in that proof.)  For the convenience of the reader, we give the following proof, which is similar to the proof of Proposition 1 in \cite[p. 87]{brm}.  Let $\mu_0 = l$ be the meromorphism on the function field $\textsf{K}=\overline{\mathbb{F}}_l(x,y)$ of the curve $E_5(b)$ induced by $[l]$; then $\mu_0$ maps $\textsf{K}$ isomorphically onto its subfield $\textsf{K}^{\mu_0}$, and $(x,y)^{\mu_0} = (f(x,y), g(x,y))$ for some rational functions $f$ and $g$ for which
\begin{equation}
[l]P = \mu_0 P = (f(P), g(P)).
\label{eq:8}
\end{equation}
If $\mathfrak{p}$ is a prime divisor of the field $\textsf{K}$, then by Hasse's formula \cite[pp. 72-73]{h2}, \cite[p. 205]{d} we have
$$[l]\mathfrak{p} = \mu_0 \mathfrak{p} = N_0(\mathfrak{p})^{\mu_0^{-1}},$$
where $N_0: \textsf{K} \rightarrow \textsf{K}^{\mu_0}$ is the norm map from $\textsf{K}$ to $\textsf{K}^{\mu_0}$.
Now the reduced norm satisfies $N(l) = l^2$ in the quaternion algebra $\textrm{End}(E_5(b))$ (see \cite[p. 199]{h}), so the defining formula $N(l) = [\textsf{K}: \textsf{K}^{\mu_0}]$ implies that the degree of inseparability of $\textsf{K}/\textsf{K}^{\mu_0}$ in characteristic $l$ is $l^2$.  In particular, the norm map $N_0$ satisfies $N_0(a) = a^{l^2}$, for elements $a \in \textsf{K}$, and similarly for divisors in $\textsf{K}$.  Now we have the divisor equality
$$(x) = \frac{\mathfrak{p}_1 \mathfrak{p}_4}{\mathfrak{o}^2},$$
where the prime divisors $\mathfrak{p}_i$ and $\mathfrak{o}$ correspond, respectively, to the points $P_i = i(0,0)$ and $O$.  Since $[l](0,0) = \pm(0,0)$, the map $[l]$ either fixes both points $P_1$ and $P_4$ or interchanges them.  Hence $\mathfrak{p}_i^{\mu_0} = N_0(\mathfrak{p}_{\pm i}) = \mathfrak{p}_{\pm i}^{l^2}$ (reading subscripts mod $5$), and therefore $x^{\mu_0} = ax^{l^2}$, for some constant $a \in \overline{\mathbb{F}}_l$ (because $\mathfrak{o}^{\mu_0} = \mathfrak{o}^{l^2}$).  Now consider $P_2 = 2(0,0) = (-b, b^2)$ on $E_5(b)$.  A similar argument shows that the divisor
$$(x+b) = \frac{\mathfrak{p}_2 \mathfrak{p}_3}{\mathfrak{o}^2}$$
satisfies $(x+b)^{\mu_0} = (x+b)^{l^2}$, implying that the field element $(x+b)^{\mu_0} = c(x+b)^{l^2}$ for some constant $c$.  Since $\mu_0$ is an isomorphism fixing constants, we have that
$$c(x+b)^{l^2} = (x+b)^{\mu_0} = x^{\mu_0}+b = ax^{l^2} + b,$$
whence it follows that 
$$c = a, \ \ \textrm{and} \ \ cb^{l^2} = b.$$
The assumption that $l \equiv 1, 4$ mod $5$ implies that $b \in \mathbb{F}_{l^2}$; hence $b^{l^2} = b$ and $a = c = 1$.  Thus $x^{\mu_0}=x^{l^2} = f(x,y)$, which implies that the endomorphism $[l] = \mu_0$ raises $X$-coordinates of points $P \in E_5(b)$ to the power $l^2$, by (\ref{eq:8}).  Hence $[l]P = \pm P^{l^2}$. \medskip

\noindent (b) Now assume that $l \equiv 3$ (mod $5$).  If the meromorphism $\mu_0$ corresponds to $[l]$, as in part (a), the proof of Theorem 6.2 in \cite[pp. 267-268]{mor} gives the following formulas:
\begin{align*}
x^{\mu_0} & = b^2 (x^{l^2}+b^{l^2}) = b^2x^{l^2}-b,\\
y^{\mu_0} & = b^3(b^{l^2}x^{l^2}+y^{l^2}) = -b^2 x^{l^2} +b^3 y^{l^2},
\end{align*}
using that $b^{l^2+1} = -1$.  A straightforward calculation yields
$$x^{\mu_0^2} = x^{l^4}, \ \ y^{\mu_0^2} = -(1+b)x^{l^4} -b - y^{l^4},$$
and this implies that
$$(x,y)^{\mu_0^2} = (x^{l^4}, -(1+b)x^{l^4} -b - y^{l^4}) = -(x^{l^4}, y^{l^4}).$$
A similar argument works if $l \equiv 2$ (mod $5$), using
\begin{align*}
x^{\mu_0} & = b^2 (x^{l^2}+b^{l^2}) = b^2x^{l^2}-b,\\
y^{\mu_0} & = -b^3(x^{l^2}+y^{l^2}+b^{l^2}) = -b^3 x^{l^2} -b^3 y^{l^2} + b^2.
\end{align*}
$\square$  \medskip

Suppose that  $l \equiv 1, 4$ mod $5$ and $g(x)$ factors into irreducible quadratics over $\mathbb{F}_l$ of the form given in Theorem 1.1B or C, with $r=\varepsilon^5(s-1)$.  In this case $\lambda \in \mathbb{F}_{l^2}$, so $\psi$ is certainly defined over $\mathbb{F}_{l^2}$.  An easy calculation shows that $\tau$ permutes the roots of $k(x) = x^2+rx+s$ and therefore $\tau(b)=b^l$ over $\mathbb{F}_l$, since $\tau(b) \neq b$ are conjugates over $\mathbb{F}_l$:
\begin{equation}
(\varepsilon^5 b+1)^2 k(\tau(b)) = 5 \sqrt{5} \varepsilon^5 k(b), \ \ \ \tau(b)=\frac{-b+\varepsilon^5}{\varepsilon^5b+1}.
\label{eq:9}
\end{equation}
(The only exception to this is $k(x)=x^2+(11 + 5\sqrt{5})x-1$, whose roots are the fixed points of $\tau(x)$.  The square $k(x)^2$ has the form of the polynomial $g(x)$ when $l \equiv \pm 1$ mod $5$; but $k(x) \mid x^4 + 22x^3 - 6x^2 - 22x + 1$, which does not have the form $g(x)$ when $l >5$.)  Then $\mu=\phi_b^l:(X,Y) \rightarrow (\phi_b(X)^l,\phi_b(Y)^l)$ is an endomorphism of $E_5(b)$, and if $P=(x,y)$,
$$\mu^2(P)=[\phi_b(\phi_b(P)^l)]^l = (\pm \phi_{b^l} \circ \phi_b(P))^{l^2} = (\pm \phi_{\tau(b)} \circ \phi_b(P))^{l^2}.$$
By Lemma 3.2 and the above remark, $\phi_{\tau(b)} \circ \phi_b = \pm [5]$, so that $\mu^2 = \pm 5l$ in $\textrm{End}(E_5(b))$, by Lemma 3.3(a).  But $\textrm{End}(E_5(b))$ is a definite quaternion algebra, so that $\mu = \phi_b^l$ satisfies $\mu^2 = -5l$ and $\mu = \pm \sqrt{-5l}$.  The same conclusion holds if $r=\bar \varepsilon^5(s-1)$ and $\bar \tau(b)=b^l$.\medskip

On the other hand, suppose that $l \equiv 2, 3$ (mod $5$), and assume $g(x)$ is an irreducible factor of $\hat H_{5,l}(X)=0$ over $\mathbb{F}_l$.  In this case the isogeny $\phi_b$ is defined over $\mathbb{F}_l(b)=\mathbb{F}_{l^4}$, where $b$ is a root of $g(x)=0$.  The map $\sigma=\left(b \rightarrow \tau(b)= \frac{-b+\varepsilon^5}{\varepsilon^5 b+1}, \sqrt{5} \rightarrow -\sqrt{5}\right)$ is an automorphism of order $4$ of $\mathbb{F}_{l^4}/\mathbb{F}_l$, satisfying $\sigma^2=(b \rightarrow -1/b)$.  Hence, $\tau(b)=b^l$ or $\sigma^3(b) = -1/\tau(b) = \bar \tau(b) = b^l$.  \medskip

In the first case, the isogeny $ \bar \phi_b: E_5(b) \rightarrow E_5(\bar \tau(b)) = E_5(b^{l^3})$ is obtained by replacing $\sqrt{5}$ by $-\sqrt{5}$ in the formulas for $\phi_b$.  The isogeny $\bar \phi_{\bar \tau(b)}: E_5(b^{l^3}) \rightarrow E_5(b)$ can be obtained by replacing $b$ by $\bar \tau(b)$ in the formulas for $\bar \phi_b$.  If the coefficients $c$ in $\bar \phi_b$ are replaced by $c^{l^3}$, this doesn't result in $\bar \phi_{\bar \tau(b)}$, but in the isogeny $\phi_{\bar \tau(b)}$ taking $E_5(\bar \tau(b))$ to $E_5(\tau(\bar \tau(b)))=E_5(-1/b)=E_5(b^{l^2})$.  We have that $\textrm{ker}(\phi_{\bar \tau(b)}) = \langle (0,0) \rangle$ in $E_5(\bar \tau(b))$, as before.  Hence, if $\mu=\bar \phi_b^l$, we have that $\mu \in \textrm{End}(E_5(b))$, and the endomorphism $\mu^2$ given by
$$\mu^2(P)=[\bar \phi_b(\bar \phi_b(P)^l)]^l = (\pm \phi_{\bar \tau(b)} \circ \bar \phi_b(P))^{l^2}$$
still has the kernel $E_5(b)[5]$.  Iterating the last formula and noting that $\sqrt{5}^{l^2} = \sqrt{5}$ gives
\begin{align*}
\mu^4(P) & = (\phi_{\bar \tau(b)} \circ \bar \phi_b(\phi_{\bar \tau(b)} \circ \bar \phi_b(P))^{l^2}))^{l^2}\\
& = ((\phi_{\bar \tau(b^{l^2})} \circ \bar \phi_{b^{l^2}}) \circ (\phi_{\bar \tau(b)} \circ \bar \phi_b)(P))^{l^4}\\
& = (\psi_2 \circ \psi_1(P))^{l^4}.
\end{align*}
Since the isogeny $\psi_1 = \phi_{\bar \tau(b)} \circ \bar \phi_b$ from $E_5(b)$ to $E_5(b^{l^2})$ has degree $5^2$ and kernel $E_5(b)[5]$, it follows that the image of $E_5(b)[25]$ under this isogeny is $E_5(b^{l^2})[5]$, which is the kernel of the degree-$5^2$ map $\psi_2 = \phi_{\bar \tau(b^{l^2})} \circ \bar \phi_{b^{l^2}}$.  (Note that $\bar \phi_b(\bar P) = \pm (0,0)$ for the point $\bar P \in E_5(b)[5]$ obtained by replacing $\sqrt{5}$ by $-\sqrt{5}$ in the formulas for the point $P$ in the proof of Lemma 3.2.)  Hence, the kernel of $\psi_2 \circ \psi_1$ is $E_5(b)[25]$.  It follows as before, on appealing to Lemma 3.3(b), that $\mu^4 = \pm 5^2 l^2$.  If $\mu^4 = -5^2 l^2$, then $\mu \in \textrm{End}(E_5(b))$ is the root of the irreducible quartic polynomial $X^4+5^2 l^2$ (irreducible over $\mathbb{Q}$ by Capelli's theorem \cite[p. 87]{co1}), which is impossible.  Therefore, $\mu^4 = 5^2 l^2$, so $\mu^2 = \pm 5l$, and we conclude that $\mu^2 = -5l$ in $\textrm{End}(E_5(b))$ (in characteristic $l > 7$).  A similar argument works if $\bar \tau(b) = b^l$.  This proves the following.  \bigskip

\noindent {\bf Theorem 3.4.} {\it If $l>7$ is a prime, then for any root $b$ of an irreducible factor of $\hat H_{5,l}(x)$ over $\mathbb{F}_l$ of the form
$$g(x) = x^4 +ax^3 +(11a+2)x^2 -ax+1,$$
or of the form $k(x)=x^2+rx+s$ with $r=\varepsilon^5(s-1)$ or $r = \bar \varepsilon^5(s-1)$, there is an endomorphism (multiplier) $\mu \in \textrm{End}(E_5(b))$ satisfying $\mu^2 = -5l$.} \bigskip

Now Deuring's lifting theorem \cite{d} yields the following theorem. \bigskip

\noindent {\bf Theorem 3.5.} {\it If $l>7$ is a prime, any irreducible factor of $\hat H_{5,l}(x)$ over $\mathbb{F}_l$ of the form
$$g(x) = x^4 +ax^3 +(11a+2)x^2 -ax+1 \ \  \textrm{or} \ \ k(x) =x^2+rx+s,$$
with $r=\varepsilon^5(s-1)$ or $r = \bar \varepsilon^5(s-1)$, arises by reduction as a factor of
$$F_{d}(x)=x^{5\textsf{h}(-d)}(1-11x-x^2)^{\textsf{h}(-d)}H_{-d}(j(x)) \ \ \textrm{modulo} \ l,$$
with
$$j(x)=\frac{(x^4+12x^3+14x^2-12x+1)^3}{x^5(1-11x-x^2)},$$
for one of the discriminants $-d=-5l$ or $-d=-20l$.} \medskip

\noindent {\it Proof.}  Let $b$ be a root of $g(x)$, resp. $k(x)$, in $\mathbb{F}_{l^4}$, resp. $\mathbb{F}_{l^2}$.  The result of Theorem 3.4 says that in characteristic $l$ there is an meromorphism (multiplier) $\bar \mu$ of the function field $\overline{\textsf{K}}$ of the curve $E_5(b)$ satisfying $\bar \mu^2 = -5l$.  Deuring's theorem \cite[p. 259]{d} says that $\overline{\textsf{K}}$ and $\bar \mu$ arise from an elliptic curve $E$ in characteristic 0 and a multiplier $\mu$ of its function field $\textsf{K}$ by reduction modulo a prime divisor of $l$:
$$\textsf{K} \rightarrow \overline{\textsf{K}}, \ \  \mu \rightarrow \bar \mu.$$
Then the characteristic-zero meromorphism ring $\textsf{M}(\textsf{K}) \cong \textrm{End}(E)$ must contain the element $\mu \notin \mathbb{Z}$.  Since algebraic relations are preserved by the injection of $\textsf{M}(\textsf{K}) \rightarrow \textsf{M}(\overline{\textsf{K}})$, it follows that $\sqrt{-5l} \in \textsf{M}(\textsf{K})$.  Hence, the curve $E$ has complex multiplication by one of the orders $\mathcal{O} = \mathbb{Z}(\sqrt{-5l})$ or $\mathbb{Z}(\frac{1+\sqrt{-5l}}{2})$.  By the theory of complex multiplication, the $j$-invariant $j(E)$ must then be a root of $H_{-5l}(X)$ or $H_{-20l}(X)$; since $j(E)$ reduces to $\bar j = j(E_5(b)) = j(b)$ in characteristic $l$, $\bar j$ is a root of $H_{-5l}(X)$ or $H_{-20l}(X)$ modulo $l$.  This proves the theorem.  $\square$ \bigskip 

The factors $k(x)=x^2+(11 \pm 5\sqrt{5})x-1$ mentioned above are also covered by the statements in these two theorems, since
\begin{align*}
F_{20}(x) & = \ (x^4 + 22x^3 - 6x^2 - 22x + 1)(x^{20} + 50x^{19} + 1150x^{18} + 14550x^{17}\\
& + 118525x^{16} + 1746272x^{15} + 34835400x^{14} + 376573200x^{13}\\
& + 1950875650x^{12} + 4311023700x^{11} + 2400976244x^{10} - 4311023700x^9\\
& + 1950875650x^8 - 376573200x^7 + 34835400x^6 - 1746272x^5 \\
& + 118525x^4 - 14550x^3 + 1150x^2 - 50x + 1),
\end{align*}
when $\epsilon_{20}=1$, i.e. when $l \equiv \pm 1$ (mod $5$) and $l \equiv 3$ (mod $4$). \medskip

\section{Analyzing the special factors of $K_{5l}(X)$}

By Theorem 3.5, the $j$-invariants corresponding to factors of the form $g(x)$ or $k(x)$ of $\hat H_{5,l}(x)$ over $\mathbb{F}_l$ are reductions of roots of $K_{5l}(X)$.  Hence, by Theorem 2.1, either
\begin{equation}
g(x) \ \textrm{or} \ k(x) \mid F_d(x) \ (\textrm{mod} \ l), \ \ \textrm{for some} \ d \in \{20\} \cup \mathfrak{T} \ \textrm{with} \ \epsilon_d = 1,
\label{eq:10}
\end{equation}
or
\begin{equation}
g(x) \ \textrm{or} \ k(x) \mid x^{10}(1-11x-x^2)^2 (j(x)^2+a_ij(x)+b_i) \ (\textrm{mod} \ l);
\label{eq:11}
\end{equation}
where $F_d(x)$ and $j(x)$ are given in Theorem 3.5 and $(X^2+a_i X+b_i)^2$ is one of the factors of $K_{5l}(X)$ in Theorem 2.1. \medskip

In this section we will show that each of the factors $H_{-d}(X)^{4\epsilon_d}$, for $d \in \mathfrak{T} - \{4\}$ in Theorem 2.1, yields $\textrm{deg}(H_{-d}(X))$ factors of the form $g(x)$, when $l \equiv 2,3$ mod $5$ (Proposition 4.1), and $2 \cdot \textrm{deg}(H_{-d}(X))$ factors of the form $k(x)$ when $l \equiv 1, 4$ mod $5$ (Proposition 4.3); always under the assumption that $\epsilon_d = 1$.  Otherwise, $\epsilon_d = 0$ and the factor corresponding to $d$ obviously contributes no factors.\medskip

\subsection{Primes $l \equiv 2, 3$ mod $5$.}

We start with the case $l \equiv 2, 3$ mod $5$. \bigskip

\noindent {\bf Proposition 4.1.} {\it When $l \equiv 2, 3$ (mod $5$) and $l \in \mathcal{S}$ or $l >379$, each of the factors $H_{-d}(X)^{4\epsilon_d}$ of $K_{5l}(X)$ in Theorem 2.1 contributes exactly $\epsilon_d \textrm{deg}(H_{-d}(X))$ irreducible factors of the form $g(x)$ in Theorem 3.5, by means of (\ref{eq:10}), to the factorization of $\hat H_{5,l}(x)$ (mod $l$), except for the factor $H_{-4}(X)^{4\epsilon_4}$ when $l \equiv 3$ (mod $4$).} \medskip

\noindent {\it Proof.}  We determine when (\ref{eq:10}) holds for the factors $H_{-d}(X)$ and $d \in \mathfrak{T}$.  First, $\epsilon_{20}=0$ for the prime $l$, so we can ignore $F_{20}(x)$ for these primes.  For the linear class equations, we have
\begin{align*}
F_4(x) = & \ (x^2 + 1)^2(x^4 + 18x^3 + 74x^2 - 18x + 1)^2,\\
F_{11}(x) = & \ (x^4 + 4x^3 + 46x^2 - 4x + 1)\\
& \times (x^8 + 32x^7 + 300x^6 + 32x^5 - 8026x^4- 32x^3 + 300x^2 - 32x + 1),\\
F_{16}(x) = & \ (x^4 + 18x^3 + 200x^2 - 18x + 1)\\
& \times (x^8 + 18x^7 - 50x^6 + 18x^5 + 15774x^4 - 18x^3 - 50x^2 - 18x + 1),\\
F_{19}(x) = & \ (x^4 + 36x^3 + 398x^2 - 36x + 1)(x^8 + 76x^6 - 24474x^4 + 76x^2 + 1).
\end{align*}
The factor $(x^2+1)^2$ is a reducible polynomial of the form $g(x)$, while the last three quartics are irreducible polynomials of the form $g(x)$ in characteristic zero.  These polynomials are denoted by $Q_d(x)$ in \cite[Prop. 4.1]{mor2}, and are irreducible factors of $F_d(x)$ of degree $4\textsf{h}(-d)$, for any $d>4$ satisfying $d \equiv \pm 1$ (mod $5$).

The discriminants of the three quartics $Q_d(x)$, for $d=11, 16, 19$ (listed in the first two rows of Table 4), are not divisible by any prime greater than $19$.  If  $\epsilon_d=1$ for a prime $l >19$, then their reductions (mod $l$) must divide $\hat H_{5,l}(x)$.  But \cite[Thm. 6.2]{mor} asserts that $\hat H_{5,l}(x)$ can only have irreducible factors which are either quartic or $x^2+1$.  If one of these factors were reducible (mod $l$), then it would have to factor as $(x^2+1)^2$ (mod $l$), which is impossible for $l >19$.  Hence, each of these quartics yields an irreducible $g(x)$ mod $l$ for any prime $l$ for which the corresponding $\epsilon_d=1$. \medskip

On the other hand, the remaining polynomials (of degrees $4$ or $8$) in the above factorizations are not divisible by polynomials over $\mathbb{F}_l$ of the form $g(x)$ for large enough primes, by the following argument.  This is obvious for the quartic factor of $F_4(x)$, when $l$ does not divide $11 \cdot 18 + 2 - 74 =2 \cdot 3^2 \cdot 7$. 
If
$$f_{11}(x) = x^8 + 32x^7 + 300x^6 + 32x^5 - 8026x^4- 32x^3 + 300x^2 - 32x + 1$$
is divisible (mod $l$) by a factor of the form $g(x)=x^4+tx^3+(11t+2)x^2-tx+1$, then since $x^8 f_{11}(-1/x)= f_{11}(x)$ and $x^4 g(-1/x) = g(x)$, we can write
\begin{align}
\notag  f_{11}(x) = & \ x^4 \tilde f_{11}\left(x-\frac{1}{x}\right), \ \ \tilde f_{11}(x) = x^4 + 32x^3 + 304x^2 + 128x - 7424,\\
\label{eq:12} g(x) = & \ x^2 \tilde g\left(x-\frac{1}{x}\right), \ \ \ \tilde g(x) = x^2 + t x + 11t + 4.
\end{align}
Then $g(x) \mid f_{11}(x)$ implies $\tilde g(x) \mid \tilde f_{11}(x)$ (mod $l$).  On the other hand, the remainder on dividing $\tilde f_{11}(x)$ by $\tilde g(x)$ is
$$r_{11}(x) = (-t^3 + 54t^2 - 648t)x - 11t^3 + 469t^2 - 3128t - 8624.$$
If this is $0$ (mod $l$), then $l$ must divide the resultant of the coefficients in $t$, which is
$$\textrm{Res}_t(-t^3 + 54t^2 - 648t, - 11t^3 + 469t^2 - 3128t - 8624) = -2^{17} 7^3 \cdot 11 \cdot 13 \cdot 19 \cdot 43.$$
It follows that $f_{11}(x)$ contributes no irreducible factors of the form $g(x)$ to $\hat H_{5,l}(x)$, for $l > 43$.  The same argument works for the $8$-th degree factors $f_{16}(x)$ and $f_{19}(x)$ in the above factorizations, since the corresponding resultants are
\begin{align*}
\textrm{Res}_t&(-t^3 + 40t^2 - 144t, -11t^3 + 315t^2 + 666t + 15876) = 2^6 3^9 7^3 \cdot 19 \cdot 43 \cdot 67;\\
\textrm{Res}_t&(-t^3 + 22t^2 - 72t, -11t^3 + 117t^2 - 792t - 24624) = 2^{17} 3^9 \cdot 13 \cdot 19 \cdot 67.
\end{align*}
\medskip
This proves that there is one irreducible factor of the form $g(x)$ dividing $\hat H_{5,l}(x)$ for each of the factors $H_{-d}(X)^{4\epsilon_d}$ (when $\epsilon_d=1$), for $d = 11, 16, 19$, and no such factors for $H_{-4}(X)^{4\epsilon_4}$, for $l > 67$.  For $l \equiv 1$ (mod $4$), this gives one factor $g(x)$ for each linear factor $H_{-d}(x)$ dividing $K_{5l}(x)$; and for $l \equiv 3$ (mod $4$), one less than that, since $H_{4}(x)^4$ divides $K_{5l}(x)$ but does not yield a factor $g(x)$.  \medskip

Consider next the quadratic class equations $H_{-d}(X)$, for $d$ one of the integers in $\{24,36,51,64,91,99\}$.  We see first that
\begin{align*}
F_{24}(x) & =  \ (x^8 - 12x^7 + 16x^6 + 3156x^5 + 16878x^4 - 3156x^3 + 16x^2 + 12x + 1)\\
& \times (x^{16} + 84x^{15} + 3236x^{14} + 73860x^{13} + 983188x^{12} + 6801300x^{11}\\
& + 18487964x^{10} + 6727524x^9 + 78893398x^8 - 6727524x^7 + 18487964x^6\\
&  - 6801300x^5 + 983188x^4 - 73860x^3 + 3236x^2 - 84x + 1)\\
& = Q_{24}(x) f_{24}(x),
\end{align*}
where $\epsilon_{24}=1$ when $\left(\frac{-6}{l}\right)=\left(\frac{2}{l}\right) = -1$, so that $\left(\frac{-3}{l}\right) = +1$.  Now the $8$-th degree factor $Q_{24}(x)$ is the product of
$$g_{24}(x) = x^4+(-6+6\sqrt{-3})x^3+(-64+66\sqrt{-3})x^2+(6-6\sqrt{-3})x+1$$
and its conjugate over $\mathbb{Q}$.  Since $\left(\frac{-3}{l}\right) = +1$, these two factors give two irreducible factors of the form $g(x)$ mod $l$, when $l > 19$, the largest prime dividing $\textrm{disc}(Q_{24}(x))$.  On the other hand, transforming $f_{24}(x)$ as we did above to yield the $8$-th degree polynomial
\begin{align*}
\tilde f_{24}(x) = & \ x^8 + 84x^7 + 3244x^6 + 74448x^5 + 1002624x^4 + 7171776x^3\\
& \ + 22449856x^2 + 27501312x + 117842176,
\end{align*}
the remainder of $\tilde f_{24}(x)$ on dividing by $\tilde g(x)$ is $A_{24}(t)x+B_{24}(t)$, with
\begin{align*}
A_{24}(t) = & \ -t^7 + 150t^6 - 9050t^5 + 280932t^4 - 4646880t^3 + 36376128t^2\\
& \ - 86966784t,\\
B_{24}(t) = & \ - 11t^7 + 1525t^6 - 83550t^5 + 2295377t^4 - 31830180t^3 + 178683408t^2\\
& \  - 134106624t + 43877376.
\end{align*}
We compute that
$$\textrm{Res}_t(A_{24}(t),B_{24}(t)) = -2^{47}3^{23}13^4 17 \cdot 19^5 23^3 37 \cdot 41 \cdot 43^2 47 \cdot 61 \cdot 67 \cdot 71 \cdot 89 \cdot 109 \cdot 113.$$
Hence, $F_{24}(x)$ contributes exactly $2 =\textrm{deg}(H_{-24}(x))$ irreducible quartics of the form $g(x)$, for $l > 113$, when $\epsilon_{24}=1$. \medskip 

The same pattern of argument works for the remaining integers in the set $\{24,36,51,64,91,99\}$.  For each such $d$, $F_d(x)$ contributes exactly $2 =\textrm{deg}(H_{-d}(x))$ irreducible quartics of the form $g(x)$, for primes $l$ greater than the largest of the primes dividing the corresponding resultants $R(d) = \textrm{Res}(A_d(t),B_d(t))$, where $r_d=A_d(t)x + B_d(t)$ is the remainder on dividing the corresponding $8$-th degree polynomial $\tilde f_d(x)$ by the polynomial $\tilde g(x)$ from (\ref{eq:12}).  The results are given in Tables 2 and 3.  Table 2 contains the quartic polynomials which are factors of each $Q_d(x)$ and reduce to polynomials of the form $g(x)$ mod $l$, and Table 3 lists the factorizations of the resultants $R(d)$.  Table 4 contains the discriminants of the polynomials $Q_d(x)$.
\medskip

\begin{table}
  \centering 
  \caption{Quartic factors $g_d(x)$ over $\mathbb{F}_l$.}\label{tab:1}

\noindent \begin{tabular}{|c|l|c|}
\hline
   &   \\
$d$	&   $\ \ \ g_d(x)=x^4+a x^3 + (11a+2)x^2 - ax+1$  \\
\hline
& \\
11 & $x^4 + 4x^3 + 46x^2 - 4x + 1$\\
& \\
16 & $x^4 + 18x^3 + 200x^2 - 18x + 1$ \\
& \\
19 & $x^4 + 36x^3 + 398x^2 - 36x + 1$ \\
\hline
  &    \\
 24  &  $x^4-(6-6\sqrt{-3})x^3+(-64+66\sqrt{-3})x^2+(6-6\sqrt{-3})x+1$  \\
  &  \\
 36  &  $x^4+(30+22\sqrt{-3})x^3+(332+242\sqrt{-3})x^2-(30+22\sqrt{-3})x+1$  \\
  &  \\
 51 & $x^4-(12-48\sqrt{-3})x^3+(-130+528\sqrt{-3})x^2+(12-48\sqrt{-3})x+1$ \\
  &  \\
 64 & $x^4-(108-63\sqrt{-2})x^3+(-1186+693\sqrt{-2})x^2+(108-63\sqrt{-2})x+1$ \\
  &  \\
 91 & $x^4-(108-144\sqrt{-7})x^3+(-1186+1584\sqrt{-7})x^2+(108-144\sqrt{-7})x+1$ \\
  & \\
 99 & $x^4+(436+176\sqrt{-3})x^3+(4798+1936\sqrt{-3})x^2-(436+176\sqrt{-3})x+1$ \\
  \hline
  &  \\
  84 & $x^4-(117+57\sqrt{-3} + 27 \sqrt{21} +33\sqrt{-7})x^3-(1285+627\sqrt{-3}$\\
  & $+297\sqrt{21}+363\sqrt{-7})x^2+(117+57\sqrt{-3} + 27 \sqrt{21} +33\sqrt{-7})x+1$\\
  &  \\
  96 & $x^4+(81+159\sqrt{-2}+129\sqrt{-3}+33\sqrt{6})x^3 + (893+1749\sqrt{-2}$\\
  &  $+1419\sqrt{-3}+363\sqrt{6})x^2-(81+159\sqrt{-2}+129\sqrt{-3}+33\sqrt{6})x+1$\\
  \hline
\end{tabular}
\end{table}

\begin{table}
  \centering 
  \caption{Prime factorization of $R(d) = \textrm{Res}(A_d(t),B_d(t))$.}\label{tab:2}

\noindent \begin{tabular}{|c|l|c|}
\hline
   &   \\
$d$	&   $\ \ \ R(d)$  (Primes with $\epsilon_d=1$ in bold.)\\
\hline
 & \\
20 & $2^{69}5^{27} {\bf 11}^9 13^8 17^4 {\bf 19}^6 {\bf 31}^4 37^3 \cdot 53 \cdot {\bf 59} \cdot {\bf 71} \cdot 73^2 \cdot {\bf 79}^2 \cdot 97$\\
  &    \\
 24  &  $-2^{47}3^{23}{\bf 13}^4 17 \cdot {\bf 19}^5 23^3 {\bf 37} \cdot 41 \cdot {\bf 43}^2 47 \cdot {\bf 61} \cdot {\bf 67} \cdot 71 \cdot 89 \cdot {\bf 109} \cdot 113$  \\
  &  \\
 36  &  $2^{52} 3^{10} {\bf 7}^{11} 11^6  {\bf 19}^3 23 \cdot  {\bf 31} \cdot  {\bf 43}^2 \cdot  {\bf 67} \cdot 71 \cdot  {\bf 79} \cdot 83 \cdot 107 \cdot  {\bf 127} \cdot  {\bf 139} \cdot  {\bf 151} \cdot  {\bf 163}$  \\
  &  $\cdot 167$\\
 51 & $-2^{75}3^{24} {\bf 7}^{11}17^2 \cdot  {\bf 31} \cdot  {\bf 37} \cdot 47^2 53 \cdot 59 \cdot  {\bf 61} \cdot  {\bf 79} \cdot 83 \cdot  {\bf 139} \cdot  {\bf 163} \cdot 179 \cdot  {\bf 211}$ \\
  &  \\
 64 & $2^{14}3^{34}7^{12}  {\bf 11}^3  {\bf19}^5 \cdot 23 \cdot 31^2  {\bf 43}^2 \cdot  {\bf 59} \cdot  {\bf 67} \cdot 79 \cdot 127 \cdot  {\bf 139} \cdot 151 \cdot  {\bf 163} \cdot 167 \cdot  {\bf 211}$ \\
  &  $ \cdot 223 \cdot  {\bf 283} \cdot  {\bf 307}$\\
 91 & $2^{74}3^{37}7^7  {\bf 11}^3 13^4 \cdot 17 \cdot  {\bf 37} \cdot 61^2 \cdot  {\bf 67} \cdot  {\bf 71} \cdot 103 \cdot  {\bf 109} \cdot 139 \cdot  {\bf 151} \cdot  {\bf 163} \cdot 283 \cdot  {\bf 331}$ \\
  &  $\cdot  {\bf 379}$ \\
 99 & $2^{75} {\bf 7}^{11} 11 \cdot  {\bf 13}^4 17 \cdot  {\bf 19}^3 29^3 41^2  {\bf 43}^2 \cdot  {\bf 61} \cdot  {\bf 79} \cdot 83^2 \cdot 107^2 \cdot  {\bf 109} \cdot  {\bf127} \cdot  {\bf 139} \cdot  {\bf 211} $ \\
  & $\cdot 227 \cdot  {\bf 283} \cdot  {\bf 307} \cdot 347$\\
  \hline
  &  \\
  84 & $2^{192}3^{84}7^{24}13^{22} 29^{10}  {\bf 43}^8 47^5 53 \cdot 59^5 61^3  {\bf 67}^3 73^2  {\bf 79}^2  83^5 97^3 \cdot 113 \cdot  {\bf 127} \cdot 131^2 \cdot  {\bf 137} $\\
  &   $\cdot 149 \cdot {\bf 151} \cdot 157 \cdot {\bf 163} \cdot 167^2 \cdot 181\cdot 197 \cdot {\bf 211}^2 \cdot 227 \cdot 229 \cdot 233^2 \cdot 241 \cdot 281 \cdot 311$\\
  & $\cdot 313 \cdot 317 \cdot {\bf 331} \cdot 349 \cdot {\bf 379} \cdot 383 \cdot 389 \cdot 397 \cdot 401 \cdot 409$\\
  &  \\
  96 & $-2^{104}3^{89} 13^{22} 17^6 {\bf 19}^{14} 23^{12} 37^7 41^6 {\bf 43}^8 47^4 61^6 {\bf 67}^7 71^2 \cdot 89 \cdot 109^3 113^2 137^2 {\bf 139}^4 \cdot {\bf 163}$\\
  & $\cdot 167 \cdot {\bf 211} \cdot 229 \cdot 239\cdot 257 \cdot 263^2 \cdot 277^2 \cdot 281 \cdot {\bf 283} \cdot {\bf 307} \cdot {\bf 331} \cdot 349 \cdot 353  \cdot 359 $\\
  & $\cdot 373 \cdot {\bf 379} \cdot 383 \cdot 397 \cdot 401 \cdot 421 \cdot 431 \cdot 449$\\
  \hline
\end{tabular}
\end{table}

\begin{table}
  \centering 
  \caption{Discriminant of $Q_d(x)$.}\label{tab:3}
  
\noindent \begin{tabular}{|c|l|c|c|}
\hline
   &  &  &  \\
$d$	&   $\textrm{disc}(Q_d(x))$  & $d$ & $\textrm{disc}(Q_d(x))$\\
\hline
 &  &  &  \\
 11 & $2^{12}5^3 11^2$ & 19 & $2^{12} 3^4 5^3 19^2$\\
 &  &  &  \\
 16 & $2^6 3^4 5^3 11^4$ &  & \\
 \hline
  &    &  &  \\
 24  &  $2^{44} 3^{16} 5^{18} 19^4$  &  64 & $2^{18}3^{24}5^{18}7^{8}11^8 19^4 59^4$ \\
  &    &  &  \\
 36  &  $2^{48} 3^{6}5^{18} 7^4 11^8 31^4$  & 91 & $2^{64}3^{24}5^{18}7^4 11^8 13^4 71^4$\\
  &    &  &  \\
 51 & $2^{64}3^{16} 5^{18} 7^{4} 31^4$ & 99 & $2^{64} 3^4 5^{18} 7^{4} 11^{12} 19^8 79^4$ \\
  \hline
  &    &  &  \\
 84 & $2^{192} 3^{64} 5^{84} 7^{20} 13^{16} 29^8 59^8 79^4$ & 96 & $2^{96} 3^{72} 5^{84} 13^{16} 17^{16} 19^8 41^8 61^8 71^8$\\
  \hline
\end{tabular}
\end{table}

A similar argument works for $d=84$ and $d=96$.  The condition $\epsilon_{84}=1$ implies that
$$\left(\frac{-84}{l}\right) = \left(\frac{3}{l}\right) = \left(\frac{7}{l}\right) = -1,$$
hence $l \equiv 3$ (mod $4$) and
$$\left(\frac{21}{l}\right) = \left(\frac{-3}{l}\right) = \left(\frac{-7}{l}\right) = +1.$$
Thus, the polynomial listed in Table 2 for $d=84$ reduces to an irreducible $g(x)$ (mod $l$), and for $l > 79$, taking its conjugates (over $\mathbb{Q}$) gives four distinct quartic factors of $\hat H_{5,l}(x)$ of this form.  A similar argument applies for $d=96$, since in this case
$$\left(\frac{-96}{l}\right) = \left(\frac{2}{l}\right) = \left(\frac{3}{l}\right) = -1,$$
in order for $\epsilon_{96}=1$, which implies once again that $l \equiv 3$ (mod $4$) and
$$\left(\frac{6}{l}\right) = \left(\frac{-2}{l}\right) = \left(\frac{-3}{l}\right) = +1.$$
To complete the argument in these cases, we transform the cofactor $f_d(x)$ in the factorization $F_d(x) = Q_d(x) f_d(x)$, obtaining $\tilde f_d(x)$ of degree $16$, and divide by the polynomial $\tilde g(x)$ from (9) to obtain the remainder $r_d=A_d(t)x+B_d(t)$, as before.  Then the resultant $R(d)$ for these two cases is given in Table 3.  Moreover, we only need to exclude the primes in bold in Table 3, since the non-bold primes have $\epsilon_d=0$.  Keeping Theorem 2.1 and Corollary 2.2 in mind, this proves Proposition 4.1. $\square$  \bigskip

Looking ahead to Sections 5 and 6, if $l \equiv 2, 3$ (mod $5$) and we can prove that each of the remaining factors $(X^2+a_i X+b_i)^2$ in Theorem 2.1 contributes only one irreducible factor of the form $g(x)$, then this will prove Theorem 1.1A for $l \in \mathcal{S}$ or $l > 379$.  Namely, the above discussion gives exactly $N=\frac{1}{4}\textrm{deg}(K_{5l}(X)) = \frac{h(-5l)}{4}$ factors $g(x)$ when $l \equiv 1$ mod $4$, and $N-1 =\frac{1}{4}(\textsf{h}(-5l)+\textsf{h}(-20l)-4)$ such factors when $l \equiv 3$ mod $4$.  Recall from the introduction that $\textsf{h}(-d)$ is the class number of the order $\textsf{R}_{-d}$ of discriminant $-d$ in $K = \mathbb{Q}(\sqrt{-5l})$, while $h(-5l)$ denotes the class number of the maximal order in $K$.  For the primes $l \equiv 3$ (mod $4$), $\textsf{h}(-5l)=h(-5l)$ and $\textsf{h}(-20l) = h(-5l)$ or $3h(-5l)$ according as $\left(\frac{-5l}{2}\right) = +1$ or $-1$, which implies the formulas of Theorem 1.1A.  This yields all the irreducible factors $g(x)$ of $\hat H_{5,l}(x)$ by Theorem 3.5. \medskip

\subsection{Primes $l \equiv 1, 4$ mod $5$.}

We turn now to the primes $l \equiv 1, 4$ (mod $5$). \bigskip

\noindent {\bf Proposition 4.2.} {\it If $l \equiv \pm 1$ (mod $5$) and $l > 79$, then for $d \in \mathfrak{T} - \{4\}$, the polynomial $g_d(x)$ in Table 2 is a product of two irreducible polynomials over $\mathbb{F}_l$ of the form $k(x) = x^2+r x+s$, where $r = \varepsilon^5(s-1)$ or $r = \bar \varepsilon^5(s-1)$.} \medskip

\noindent {\it Proof.} This follows from the computations at the beginning of Section 3, according to which the roots of $g_d(x)=x^4+ax^3+(11a+2)x^2-ax+1$ are linear combinations of the square-roots $\sqrt{-\Theta_i}$, where $\Theta_2 \Theta_3 = -\frac{a^2}{4} \Theta_1$, and
$$-\Theta_1 = \frac{1}{4}(a^2-44a-16).$$
Each of the $\Theta_i$ lies in $\mathbb{F}_l$, for primes $l \equiv \pm 1$ (mod $5$), and none is zero in $\mathbb{F}_l$, since the largest prime factor of any value $N_\mathbb{Q}(a)$ in Table 5 is 19, and the largest prime factor of any value $N_\mathbb{Q}(a^2-44a-16)$ is $79$.  Moreover, the final column in Table 5 shows that in every case, $\displaystyle \left(\frac{-\Theta_1}{l}\right) = \left(\frac{-d}{l}\right) = -1$, assuming $\epsilon_d = 1$, so that $-\Theta_1$ is always a quadratic nonresidue (mod $l$).  Hence, one of $-\Theta_2$ and $-\Theta_3$ is a quadratic residue, and one is a quadratic nonresidue.  If $-\Theta_2$ is a quadratic residue, then by the formulas for the roots $\rho_i$ in Section 3,
$$x^2+rx+s = (x-\rho_1)(x-\rho_2) \ \textrm{and} \ x^2+r'x+s' = (x-\rho_3)(x-\rho_4)$$
are factors of $g_d(x)$ over $\mathbb{F}_l$ with $r = \bar \varepsilon^5(s-1)$ and $r' = \bar \varepsilon^5(s'-1)$; while if $-\Theta_3$ is a quadratic residue, then
$$x^2+rx+s = (x-\rho_1)(x-\rho_4) \ \ x^2+r'x+s' = (x-\rho_2)(x-\rho_3)$$
are factors of $g_d(x)$ over $\mathbb{F}_l$ with $r = \varepsilon^5(s-1)$ and $r' = \varepsilon^5(s'-1)$.  These factors are irreducible over $\mathbb{F}_l$ in every case, since the Pellet-Stickelberger-Voronoi Theorem (see \cite[Appendix]{brm}) implies that $g_d(x)$ has an even number of irreducible factors.  This is because $\textrm{disc}(g_d(x))=125a^2(a^2 - 44a - 16)^2$ is always a square in $\mathbb{F}_l$, and the PSV Theorem says that
$$\left(\frac{\textrm{disc}(g_d(x))}{l}\right) = (-1)^{4+r},$$
where $r$ is the number of irreducible factors (mod $l$).  We know that $r \ge 2$.  If $r= 4$, then all the roots $\rho_i \in \mathbb{F}_l$, implying that all the numbers $-\Theta_i$ would be squares (mod $l$), which is not the case.  Therefore, $r=2$ and we get two distinct factors of the form $k(x)=x^2+rx+s$ satisfying the required conditions, for each factor $g_d(x)$ in Table 2.  This proves the proposition. $\square$ \bigskip
 
 \begin{table}
  \centering 
  \caption{$-4\Theta_1 = a^2-44a-16$ for $a$ in $g_d(x)$.}\label{tab:4}
  
\noindent \begin{tabular}{|c|l|c|c|}
\hline
   &  &  \\
$d$	&   $a$  & $a^2-44a-16$ \\
\hline
 &  &   \\
 11 & $4$ & $-2^4 \cdot 11$\\
 16 & $18$ & $-2^2 \cdot 11^2$ \\
  19 &  $36$  &  $-2^4 \cdot 19$ \\
 \hline
  &    &    \\
 24  &  $-6+6\sqrt{-3}$ & $2^5 \left(\frac{7-3\sqrt{-3}}{2}\right)^2$ \\
  &    &   \\
 36  &  $30+22\sqrt{-3}$ &$-2^6 \left(\frac{11-\sqrt{-3}}{2}\right)^2$ \\
  &    &    \\
 51 & $-12+48\sqrt{-3}$ & $2^4 (17) (2-3\sqrt{-3})^2$\\
  &  &  \\
  64 & $-108+63\sqrt{-2}$ & $-(90+91\sqrt{-2})^2$\\
  &  & \\
  91 & $-108+144\sqrt{-7}$ &$2^4(13)(9-10\sqrt{-7})^2$ \\
  &  & \\
  99 & $436+176\sqrt{-3}$ & $-2^4(11)(23-18\sqrt{-3})^2$ \\
  \hline
  &    &    \\
 84 & $-117-57\sqrt{-3}-27\sqrt{21}-33\sqrt{-7}$ & $-(-104 + 78\sqrt{-3} - 18\sqrt{21} + 48\sqrt{-7})^2$\\
  & & \\
  96 & $81+159\sqrt{-2}+129\sqrt{-3}+33\sqrt{6}$ & $-(-221+51\sqrt{-2}+39\sqrt{-3}-93\sqrt{6})^2$ \\
  \hline
\end{tabular}
\end{table}

\noindent {\bf Proposition 4.3.} {\it When $l \equiv 1, 4$ (mod $5$) and $l \in \mathcal{S}$ or $l >379$, the following facts hold.} \smallskip

\noindent (a) {\it Each of the factors $H_{-d}(X)^{4\epsilon_d}$ of $K_{5l}(X)$ in Theorem 2.1, for $d \in \mathfrak{T} - \{4\}$, contributes exactly $2\epsilon_d \textrm{deg}(H_{-d}(X))$ irreducible factors of the form $k(x)$ to the factorization of $\hat H_{5,l}(x)$ (mod $l$), by means of (\ref{eq:10}).} \smallskip

\noindent (b) {\it When $l \equiv 3$ (mod $4$), $H_{-20}(X)^2$ (mod $l$) contributes $2\delta_{l,1}$ factors of the form $k(x)$ to the factorization of $\hat H_{5,l}(x)$ (mod $l$); where $\delta_{l,1} = 1$ if $l \equiv 1$ (mod $5$) and $\delta_{l,1} = 0$ if $l \equiv 4$ (mod $5$).} \smallskip

\noindent (c) {\it When $l \equiv 3$ (mod $4$), the factor $H_{-4}(X)^4$ in Theorem 2.1 only contributes the factor $k(x) = x^2+1$ to the factorization of $\hat H_{5,l}(x)$ (mod $l$).} \medskip

\noindent {\it Proof.} (a) Aside from the quadratic factors $k(x)=x^2+(11 \pm 5\sqrt{5})x-1$ of
$$Q_{20}(x) = x^4 + 22x^3 - 6x^2 - 22x + 1,$$
the factor $k(x)=x^2+rx+s$ always occurs with its companion $\bar k(x) = \frac{1}{s} x^2 k(-1/x) = x^2-\frac{r}{s}x+\frac{1}{s}$, such that $k(x) \bar k(x) = g(x)$ is a polynomial of the form $g(x)$.  Since the roots of each cofactor $f_d(x)$ are stable under the mapping $x \rightarrow -1/x$, the calculations in subsection 4.1 show that none of these cofactors is divisible by a polynomial of the form $k(x)$, for $l \in \mathcal{S}$ or $l > 379$.  The assertion now follows from Proposition 4.2. \medskip

\noindent (b) and (c):  When $l \equiv 3$ mod $4$, the factor $k(x) = x^2+1$ divides $F_4(x)$, and is certainly a factor of $\hat H_{5,l}(x)$.  The argument in (a) shows that the cofactor $f_4(x)$ of $(x^2+1)^2$ in $F_4(x)$ contributes no factors $k(x)$.  This proves (c).  Moreover, in this case $\epsilon_{20}=1$, and $F_{20}(x)$ is divisible by $Q_{20}(x)=(x^2+(11 + 5\sqrt{5})x-1)(x^2+(11 - 5\sqrt{5})x-1)$, which is a product of two factors of the form $k(x) = x^2+rx+s$.  Also, 
 $$\textrm{disc}(x^2+(11 + 5\sqrt{5})x-1) = 250 + 110\sqrt{5} = 2^2 \sqrt{5}^3 \left(\frac{1 + \sqrt{5}}{2}\right)^5.$$
 Thus, $q(x)=x^2+(11 + 5\sqrt{5})x-1$ factors into linear factors (mod $l$) if $\left(\frac{(5+\sqrt{5})/2}{l}\right) = +1$ and is irreducible (mod $l$) otherwise.  Now $L_2=\mathbb{Q}\left(\sqrt{\frac{5+\sqrt{5}}{2}}\right)$ is the real subfield of the field $\mathbb{Q}(\zeta_{20})$ (see \cite[p. 6]{mor8}).  It corresponds by class field theory to the congruence subgroup $H = \{\pm 1 \ \textrm{mod} \ 20\} \subset \left(\mathbb{Z}/20\mathbb{Z}\right)^\times$.  If $l \equiv 3$ (mod $4$) and $l \equiv 4$ (mod $5$), then $l \equiv 19$ (mod $20$) so the polynomial $Q_{20}(x)$ splits, and $d=20$ yields no factors of the form $k(x)$.  On the other hand, if $l \equiv 3$ (mod $4$) and $l \equiv 1$ (mod $5$), then $l \equiv 11$ (mod $20$) has order $2$ modulo $H$, implying that $Q_{20}(x)$ is a product of two irreducible polynomials (mod $l$).  Thus, $d=20$ yields two factors of the form $k(x)$, when $l \equiv 1$ (mod $5$).  It is only necessary to check that $f_{20}(x)$ is also not divisible by a polynomial of the form $g(x)$ modulo any of these primes, and this follows from the entry $R(20)$ in Table 3.  $\square$ \medskip
 
Looking ahead to Sections 5 and 6, if $l \equiv 1, 4$ (mod $5$) and we can show that each of the factors $(X^2+a_iX+b_i)^2$ in Theorem 2.1 yields two factors of the form $k(x)$, then we can complete the proof of Theorems 1.1B and C as follows.  For each of the factors $H_{-d}(X)^{4\epsilon_d}$ with $d \neq 4$, we have $2 \epsilon_d \textrm{deg}(H_{-d}(X))$ factors of the form $k(x)$ dividing $\hat H_{5,l}(x)$.  By Proposition 4.3, this gives a total of
 $$\sum_{d \in \mathfrak{T}-\{4\}}{2\epsilon_d \textrm{deg}(H_{-d}(X))} + \sum_{X^2+a_iX+b_i}{2}+ \epsilon_4(1 + 2\delta_{l,1}),$$
factors of the form $k(x)$.  But this latter expression equals
 $$\frac{1}{2}(\textrm{deg}(K_{5l}(X))-4\epsilon_{20}-4\epsilon_4)+ \epsilon_4(1 + 2\delta_{l,1})=\frac{a_l}{2}h(-5l)-3\epsilon_4+2\epsilon_4\delta_{l,1},$$
where the last equality follows from $\epsilon_4=\epsilon_{20}$, using the notation $a_l$ of Theorem 2.3.  This equals $\frac{a_l}{2}h(-5l)$ if $l \equiv 1$ (mod $4$); $\frac{a_l}{2}h(-5l)-3$ if $l \equiv 3$ (mod $4$) and $l \equiv 4$ (mod $5$); and $\frac{a_l}{2}h(-5l)-1$ if $l \equiv 3$ (mod $4$) and $l \equiv 1$ (mod $5$).  This yields the formulas of Theorem 1.1B and C, for $l \in \mathcal{S}$ or $l > 379$. \medskip

\section{Values of $r(\tau)$ over $\mathbb{Q}(\sqrt{-5l})$}

To prove Theorem 1.1 it suffices, by the final remarks in subsections 4.1 and 4.2, to prove that each factor $(X^2+a_iX+b_i)^2$ in Theorem 2.1 yields -- by means of (\ref{eq:11}) -- exactly one factor of the form $g(x)$, or two factors of the form $k(x)$, of the Hasse invariant $\hat H_{5,l}(x)$, for the appropriate congruence conditions on $l$.  In this section we show that each factor $(X^2+a_iX+b_i)^2$ yields at least one factor $g(x)$, respectively two factors $k(x)$, of $\hat H_{5,l}(x)$; then in Section 6 we show that no more than one, respectively two, such factors arise from $(X^2+a_iX+b_i)^2$. \medskip

Let $r(w)$ denote the Rogers-Ramanujan continued fraction:

\begin{align*}
r(w) &= \frac{q^{1/5}}{1+\frac{q}{1+\frac{q^2}{1+ \frac{q^3}{1+\cdots}}}}=\frac{q^{1/5}}{1+} \ \frac{q}{1+} \ \frac{q^2}{1+} \ \frac{q^3}{1+} \dots,\\
& = q^{1/5} \prod_{n \ge 1}{(1-q^n)^{(n/5)}},  \ \ q = e^{2 \pi i w}, \ \ w \in \mathbb{H}.
\end{align*}
See  \cite{anb},  \cite{ber}, \cite{du}.  We recall the following identities for $r(w)$ from \cite[(2.5), (7.3)]{du}:
\begin{align}
\label{eq:13} & j(w)  = \frac{(r^{20}-228r^{15}+494r^{10}+228r^5+1)^3}{r^5(1-11r^5-r^{10})^5}, \ \ r=r(w);\\
\label{eq:14} & r^5\left(\frac{-1}{5w}\right) = \frac{-r^5(w)+\varepsilon^5}{\varepsilon^5 r^5(w)+1} = \tau(r^5(w)).
\end{align}
Here and in the rest of this section we denote by $j(w)$ the modular $j$-function.  These identities imply that
\begin{equation*}
j\left(\frac{-1}{5w}\right) = j(5w) = \frac{(r^{20}+12r^{15}+14r^{10}-12r^5+1)^3}{r^{25}(1-11r^5-r^{10})}, \ \ r=r(w),
\end{equation*}
and therefore
\begin{equation}
j(w) = \frac{(r^{20}+12r^{15}+14r^{10}-12r^5+1)^3}{r^{25}(1-11r^5-r^{10})}, \ \ r=r(w/5).
\label{eq:15}
\end{equation}

\noindent {\bf Proposition 5.1.} (a) {\it If $l \equiv 3$ (mod $4$) and
$$\tau_1=\frac{5+\sqrt{-5l}}{10},$$
then $\rho_1=r(\tau_1)^5$ is a real algebraic number of degree $2h(-5l)$ over $\mathbb{Q}$, contained in the ray class field $\Sigma_{\wp_5}$, where $\wp_5^2 \cong 5$ in $K = \mathbb{Q}(\sqrt{-5l})$.  Moreover, $\Sigma_{\wp_5}=K(\rho_1)$.} \smallskip

(b) {\it Let $f=1$, if $l \equiv 1$ (mod $4$); and $f=2$, if $l \equiv 3$ (mod $4$).  If
$$\tau_2=\frac{\sqrt{-5l}}{5},$$
then $\rho_2=r(\tau_2)^5$ is a real algebraic number of degree $2\textsf{h}(-20l)$ over $\mathbb{Q}$, contained in the class field $\Sigma_{\wp_5}\Omega_f$.  Moreover, $\Sigma_{\wp_5}\Omega_f=K(\rho_2)$.} \smallskip

\noindent {\it Proof.} (a) We first recall the identity
$$\frac{1}{r^5(\tau)}-11-r^5(\tau) = \left(\frac{\eta(\tau)}{\eta(5\tau)}\right)^6,$$
where $\eta(\tau)$ is the Dedekind $\eta$-function.  From Schertz \cite[p. 159]{sch}, applied to the function $\left(\frac{\eta(w/5)}{\eta(w)}\right)^6 \gamma_3(w)^2$, with $w = 5\tau_1$, we have that $\lambda = \left(\frac{\eta(\tau_1)}{\eta(5\tau_1)}\right)^6 \in \Sigma$, the Hilbert class field of $K$, and
$$q(r^5(\tau_1),\lambda)=r^{10}(\tau_1)+(11+\lambda)r^5(\tau_1)-1=0.$$
The discriminant $\delta$ of the quadratic $q(x,\lambda)=x^2+(11+\lambda)x-1 \in \Sigma[x]$ is
$$\delta = (11+\lambda)^2+4 = 125+22\lambda+\lambda^2\cong \wp_5^3 \mathfrak{b}, \ \ (\wp_5, \mathfrak{b})=1,$$
using the fact that $\lambda \cong \wp_5^3$ from Deuring's treatise \cite[p. 43]{d2} and that $125=5^3 \cong \wp_5^6$. Now $\wp_5$ is unramified in $\Sigma/K$, so is not the square of a divisor in $\Sigma$.  It follows that $q(x, \lambda)$ is irreducible over $\Sigma$, so that $\rho_1$ generates a quadratic extension of $\Sigma$.  Also,
$$q_1=e^{2\pi i\tau_1} = e^{\pi i-\frac{\sqrt{5l}}{5}\pi} = -e^{-\frac{\sqrt{5l}}{5}\pi}$$
is real, so that $r(\tau_1)^5 \in \mathbb{R}$.  \medskip

Now by Thm. 15.16  in \cite{co} (a result of Cho \cite{cho}), the fact that $r^5(\tau)$ lies in the field $\textsf{F}_5$ of modular functions for $\Gamma(5)$ implies that $r^5(\tau_1) \in L_{\mathcal{O}_K,5} = \Sigma_5$, where $\Sigma_5$ is the ray class field for the conductor $\mathfrak{f}=(5)$.  On the other hand, $[\Sigma_5: \Sigma] = \frac{\varphi_K(5)}{2} = 10$, and therefore $\rho_1$ generates the unique quadratic subfield of $\Sigma_5/\Sigma$, which is $\Sigma_{\wp_5}$.  In particular, $\mathbb{Q}(\rho_1) \subseteq \Sigma_{\wp_5}^+$, the real subfield of the normal extension $\Sigma_{\wp_5}/\mathbb{Q}$.  Note that $[\Sigma_{\wp_5}^+: \mathbb{Q}] = 2h(-5l)$.  \medskip

By the identity (\ref{eq:13}), $j(\tau_1)$ is a rational function of $r(\tau_1)^5$.  Since $\tau_1$ is a basis quotient for an ideal in the maximal order $R_K$ of $K$, the quantity $j(\tau_1)$ generates a subfield of the Hilbert class field $\Sigma$ of degree $h(-5l)$ over $\mathbb{Q}$.  Hence
$$\mathbb{Q}(j(\tau_1)) \subseteq \mathbb{Q}(\rho_1) \subseteq \Sigma_{\wp_5}^+$$
implies that $\rho_1=r^5(\tau_1)$ has degree $h(-5l)$ or $2h(-5l)$ over $\mathbb{Q}$.  Since $\rho_1 \notin \Sigma^+=\mathbb{Q}(j(\tau_1))$, the latter must hold.  Finally, the fact that $\rho_1$ is real implies that $K(\rho_1) = K \mathbb{Q}(\rho_1) = K \Sigma_{\wp_5}^+ = \Sigma_{\wp_5}$.  This completes the proof.  \medskip

(b) The same arguments work for part (b), if $\Sigma, \Sigma_5$ and $\Sigma_{\wp_5}$ are replaced everywhere by $\Omega_2, \Sigma_5 \Omega_2$ and $\Sigma_{\wp_5}\Omega_2$, when $l \equiv 3$ (mod $4$).  In this case, $L_{\mathcal{O}_K,5}$ is also to be replaced by $L_{\mathcal{O}_2,5}$, where $\mathcal{O}_2=\textsf{R}_{-20l}$ is the order of discriminant $-20l$ in $K$ and $\textsf{h}(-20l)$ is its class number.  $\square$  \bigskip

In the next lemma we use the following notation from \cite[p. 1184]{mor2}.  Let
$$z(w) = r^5(w/5)-r^{-5}(w/5), \ \ w \in \mathbb{H},$$
so that
\begin{equation}
z(5\tau_k) = -11 - \lambda(\tau_k) = -11 - \left(\frac{\eta(\tau_k)}{\eta(5\tau_k)}\right)^6.
\label{eq:16}
\end{equation}
Note from the first paragraph of the above proof (and its analogue in part (b)) it is clear that $z(5\tau_i) \in \Omega_f$.  From \cite[pp. 1180, 1184]{mor2}, or by direct calculation using (\ref{eq:15}), we have
\begin{equation}
j(w) = -\frac{(z(w)^2+12z(w)+16)^3}{z(w)+11} = J(z(w)).
\label{eq:17}
\end{equation}

\noindent {\bf Lemma 5.2.} {\it Let $z_i = z(w_i)$, where $w_i = 5\tau_i$, $\tau_i$ as in Proposition 5.1.  For a given ideal $\mathfrak{a} = (a,w) \subseteq \textsf{R}_{-d}$ ($d=5l$ or $20l$) with ideal basis quotient $\tau=w/a$, where $(a,f)=1$ and $5a \mid N(w)$, there is a unique value of $\tilde z \in \Omega_f$ for which
$$J(\tilde z) = -\frac{(\tilde z^2+12 \tilde z+16)^3}{\tilde z+11} = j(w/a)$$
and $\tilde z+11 \cong \wp_5^3$, and this value is $\tilde z = z_i^{\sigma^{-1}}$, where $\sigma = \left(\frac{\Omega_f/K}{\mathfrak{a}R_K}\right)$.} \medskip

\noindent {\it Proof.}  We have $J(z_i) = j(w_i)$ from (\ref{eq:16}) and (\ref{eq:17}).  Furthermore, $\{1, w_i\}$ is a basis for the order $\mathcal{O}_f$, where the conductor $f = 1$ or $2$ as in Proposition 5.1.  Hence,
$$J(z_i^{\sigma^{-1}}) = J(z_i)^{\sigma^{-1}} = j(w_i)^{\sigma^{-1}} = j(\mathfrak{a}) = j(w/a),$$
from the theory of complex multiplication (see \cite[p. 125]{h3}, \cite[p. 35]{d2}, or \cite[p. 218]{co}).  Now assume there is another $\tilde z \in \Omega_f$ different from $z_i^{\sigma^{-1}}$ for which $J(\tilde z) = j(w/a)$.
As in \cite[Lemma 2.2]{mor2} we define
\begin{align*}
F(u,v) &= -uv\frac{J(u-11)-J(v-11)}{u-v}\\
&=u^5 v + u^4 v^2 + u^3 v^3 + u^2 v^4 + u v^5 - 30 u^4 v - 30 u^3 v^2 - 30 u^2 v^3 - 30 uv^4\\
& + 315 u^3 v + 315 u^2 v^2 + 315 u v^3 - 1300 u^2 v - 1300 u v^2 + 1575 u v - 125.
\end{align*}
Then $F(u,v) = 0$ is a curve of genus $0$, parametrized by the rational functions
\begin{align*}
u&=-\frac{125}{t(t^4+5t^3+15t^2+25t+25)}\\
v&=-\frac{t^5}{t^4+5t^3+15t^2+25t+25}.
\end{align*}
Hence, $F(\tilde z+11,z_i^{\sigma^{-1}}+11)=0$ gives that
$$\tilde z+11= \frac{-125}{t(t^4+5t^3+15t^2+25t+25)},$$
or
$$t^5+5t^4+15t^3+25t^2+25t+\frac{125}{\tilde z+11}=0,$$
for some algebraic number $t$.  Since $\tilde z+11 \cong \wp_5^3$, $t$ is an algebraic integer.  The equation for $v$ gives that
$$z_i^{\sigma^{-1}}+11=\frac{-t^5}{t^4+5t^3+15t^2+25t+25}=\frac{t^5}{\frac{125}{t(\tilde z+11)}}=t^6 \frac{(\tilde z+11)}{125}.$$
It follows that
$$t^6 =5^3 \frac{z_i^{\sigma^{-1}}+11}{\tilde z + 11} \cong 5^3 \cong \wp_5^6.$$
Therefore, $t \cong \wp_5$.  Furthermore, if $\mathfrak{p}$ is a prime divisor of $\wp_5$ in $\Omega'=\Omega_f(t)$, then $\mathfrak{p}$ divides the algebraic integer $t$, which implies in turn that
$$\mathfrak{p}^5 \mid \theta = t^5+5t^4+15t^3+25t^2+25t.$$
But then $\mathfrak{p}^5 \mid \frac{125}{\tilde z+11} \cong \wp_5^3$, which implies that $\mathfrak{p}^2 \mid \wp_5$, so that $\wp_5$ is ramified in $\Omega'/K$.  Let $e \ge 2$ denote the ramification index of $\mathfrak{p}$ over $\wp_5$.  Then $\mathfrak{p}^e || t$.  But this implies that $\mathfrak{p}^{5e} \mid \theta$, while $\mathfrak{p}^{3e}$ exactly divides $\frac{125}{\tilde z+11}=-\theta$.  This contradiction establishes that no such $\tilde z$, distinct from $z_i^{\sigma^{-1}}$, exists. $\square$ \bigskip

\noindent {\bf Remark.} Lemma 5.2 says that the quantities $z(w/a) = z(w_i)^{\sigma^{-1}}$ transform in the same way that the $j$-invariants $j(w/a) = j(w_i)^{\sigma^{-1}}$ transform under the Galois group $\textrm{Gal}(\Omega_f/K)$.  \medskip 

The following theorem holds the keys to proving the existence of the polynomials $g(x)$ and $k(x)$ in Theorem 1.1. \medskip

\noindent {\bf Theorem 5.3.} (a) {\it Let the function $j_5^*(\tau)$ be defined by
\begin{equation}
j_5^*(\tau)= \left(\frac{\eta(\tau)}{\eta(5\tau)}\right)^6+22+125\left(\frac{\eta(5\tau)}{\eta(\tau)}\right)^6.
\label{eq:18}
\end{equation}
If $L \subset \Omega_f$ is the fixed field of $\sigma_{\wp_5} = \left(\frac{\Omega_f/K}{\wp_5}\right) \in \textrm{Gal}(\Omega_f/K)$, then the minimal polynomial of $\rho_i = r^5(\tau_i) \in \Sigma_{\wp_5} \Omega_f$ over $L$ is}
$$g_L(x) = x^4+j_5^*(\tau_i) x^3 + (11j_5^*(\tau_i)+2)x^2 - j_5^*(\tau_i)x+1.$$

\noindent (b) {\it The prime divisor $\mathfrak{l}$, for which $\mathfrak{l}^2 = (l)$ in $K$, splits completely in the field $L$.  Thus, if $\mathfrak{q}$ is a prime divisor of $\mathfrak{l}$ in $L$, then $\mathfrak{q}$ has degree $1$ over $l$, and
$$g(x) := g_L(x) \ \textrm{mod} \ \mathfrak{q}$$
is a polynomial in $\mathbb{F}_l[x]$.} \smallskip

\noindent (c) {\it If $l \equiv 2, 3$ (mod 5), the polynomial $g_L(x)$ irreducible in $L_\mathfrak{q}[x]$, where $L_\mathfrak{q}$ is the completion of the field $L$ with respect to the prime $\mathfrak{q}$.  If $l \equiv 1, 4$ (mod 5), the polynomial $g_L(x)$ splits into two irreducible quadratics in $L_\mathfrak{q}[x]$.} \medskip

\noindent {\it Proof.}  In order to work with certain values of $j(\tau)$, where $\tau$ is the basis quotient for an ideal in the order $\textsf{R}_{-d} = \mathcal{O}_f$, we first note the following. \medskip

When $l \equiv 3$ mod $4$ and $f=1$, the quantities $\frac{-1}{\tau_1}=\frac{-5+\sqrt{-5l}}{2((5+l)/4)}$ and $\frac{-1}{5\tau_1}=\frac{-5+\sqrt{-5l}}{2((5^2+5l)/4)}$ are basis quotients for ideals $\mathfrak{a}=(\frac{5+l}{4},\frac{-5+\sqrt{-5l}}{2})$ with norm $N(\mathfrak{a})=\frac{5+l}{4}$ and $\wp_5\mathfrak{a}$ with norm $N(\wp_5\mathfrak{a})=\frac{5^2+5l}{4}$.  In this case, we have that
$$\wp_5\mathfrak{a} = \left(\frac{5(5+l)}{4},\frac{-5+\sqrt{-5l}}{2}\right) = \left(\frac{-5+\sqrt{-5l}}{2}\right).$$
When $l \equiv 1$ (mod $4$) or $l \equiv 3$ (mod $4$) and $f=2$, $\frac{-1}{\tau_2}=\frac{\sqrt{-5l}}{l}$ and $\frac{-1}{5\tau_2}=\frac{\sqrt{-5l}}{5l}$ are basis quotients for the ideals $\mathfrak{l}=(l, \sqrt{-5l})$ and $\wp_5 \mathfrak{l} =(5l, \sqrt{-5l}) = (\sqrt{-5l})$, respectively.  Let $\mathfrak{a}$ denote the ideal $\mathfrak{l}$ in this case, so that $\wp_5 \mathfrak{a} \sim 1$ (mod $f$).  It is clear that $\mathfrak{l} \sim \wp_5$ in $K$, so we also have $\mathfrak{a} \sim \mathfrak{l}$ when $l \equiv 3$ (mod $4$) and $f=1$.  \medskip

The identity (\ref{eq:14}) shows that the quantities $r^5(-1/(5\tau_k))$ are quadratic over $\Omega_f$ and lie in $\Sigma_{\wp_5}\Omega_f$.  (Recall that $\sqrt{5}$ lies in the genus field of $K$ and is therefore contained in the Hilbert class field.)  Since
$$j\left(\frac{-1}{\tau_k}\right) = j(\mathfrak{a}) = j(\mathcal{O}_f)^{\sigma_\mathfrak{a}^{-1}} = j(w_k)^{\sigma_\mathfrak{a}^{-1}},$$
where $\sigma_\mathfrak{a}=\left(\frac{\Omega_f/K}{\mathfrak{a}}\right)$, Lemma 5.2 implies that
$$z\left(\frac{-1}{\tau_k}\right) = z(w_k)^{\sigma_\mathfrak{a}^{-1}},$$
and therefore (\ref{eq:16}) gives that
\begin{equation*}
\lambda'=\lambda \left(\frac{-1}{5\tau_k}\right) = \lambda(\tau_k)^{\sigma_\mathfrak{a}^{-1}}.
\end{equation*}
Note that $\lambda'=\lambda\left(\frac{-1}{5\tau_k}\right)=\frac{5^3}{\lambda(\tau_k)}$ by the transformation formula $\eta(-1/\tau)=\sqrt{\frac{\tau}{i}} \eta(\tau)$.  Therefore,
\begin{equation}
\lambda'= \lambda(\tau_k)^{\sigma_\mathfrak{a}^{-1}} = \frac{5^3}{\lambda(\tau_k)}.
\label{eq:19}
\end{equation}
Thus, $\lambda'$ and $\lambda = \lambda(\tau_k)$ are conjugate values over $K$ and
\begin{align*}
(x^2+(11+& \lambda)x-1)(x^2+(11+\lambda')x-1) \\
 & = (x^2+(11+\lambda)x-1)(x^2+(11+\frac{5^3}{\lambda})x-1) \\
 & =x^4+(22+\lambda+\frac{5^3}{\lambda})x^3 +[11(22+\lambda+\frac{5^3}{\lambda})+2]x^2\\
 &    \ \ \ \ \ \ \ -(22+\lambda+\frac{5^3}{\lambda})x+1\\
 & = x^4+j_5^*(\tau_k) x^3 + (11j_5^*(\tau_k)+2)x^2 - j_5^*(\tau_k)x+1,
 \end{align*}
where the modular function $j_5^*(\tau)$ is defined by (\ref{eq:18}).  (See \cite{na} and \cite{mor8}.)  This, together with the proof of Proposition 5.1, implies that
$$g_L(x) = x^4+j_5^*(\tau_k) x^3 + (11j_5^*(\tau_k)+2)x^2 - j_5^*(\tau_k)x+1$$
divides the minimal polynomial of $r^5(\tau_k)$ over $K$.
\medskip

However, $\mathfrak{a} \sim \wp_5$ (mod $f$) implies that $\sigma_\mathfrak{a}=\sigma_{\wp_5}$ and $\sigma_\mathfrak{a}$ has order $2$.  Then the formula
$$j_5^*(\tau_k) = 22+\lambda+\frac{5^3}{\lambda} = -z(w_k)-z(w_k)^{\sigma_{\wp_5}}$$
implies that $j_5^*(\tau_k) \in L$, the fixed field of the subgroup $\langle \sigma_{\wp_5} \rangle$ of $\textrm{Gal}(\Omega_f/K)$.  Now $[\Sigma_{\wp_5}\Omega_f:L] = 4$ implies that $g_L(x)$ is the minimal polynomial of $\rho_k=r^5(\tau_k)$ over $L$.  The fact that the coefficients of $g_L(x)$ are linear expressions in $j_5^*(\tau_k)$ implies that
$$L = K(j_5^*(\tau_k)) \  \ \textrm{with} \ \ \langle \sigma_{\wp_5} \rangle = \textrm{Gal}(\Omega_f/L).$$
It follows from Artin Reciprocity that a prime ideal $\mathfrak{p}$ of $K$ splits in $L$ if and only if $\mathfrak{p} \sim 1$ or $\wp_5$ (mod $f$).  In particular, $\mathfrak{l} \sim \wp_5$ (mod $f$), so that:
$$\mathfrak{l} \ \textrm{splits in the field} \ L.$$
If $\mathfrak{q}$ is any prime divisor of $\mathfrak{l}$ in $L$, then $j_5^*(\tau_k)$ (mod $\mathfrak{q}$) lies in the prime field $\mathbb{F}_l$, so that the coefficients of $g_L(x)$ (mod $\mathfrak{q}$) also lie in $\mathbb{F}_l$.  This proves (a) and (b).  \medskip

If $l \equiv 2,3$ (mod $5$), then $\mathfrak{l}^2 = (l)$ implies that $\mathfrak{l}$ has order $4$ in the ray class group $\mathcal{C}_{\wp_5}$ modulo $\wp_5$ in $K$.  If $l \equiv \pm 1$ (mod $5$), then $\mathfrak{l}$ has order $2$ in $\mathcal{C}_{\wp_5}$.  In the first case, a prime divisor $\mathfrak{q}$ of $\mathfrak{l}$ in $L$ is inert in $\Sigma_{\wp_5}\Omega_f$, so that $g_L(x)$ is irreducible over the completion $L_\mathfrak{q}$, and therefore factors mod $\mathfrak{q}$ as a power of an irreducible polynomial over $R_L/\mathfrak{q} \cong \mathbb{Z}/l\mathbb{Z}$.  In the second case, $\mathfrak{q}$ splits into two primes of relative (and absolute) degree $2$ in $\Sigma_{\wp_5}\Omega_f$, so that $g_L(x)$ splits into two irreducible quadratics over $L_\mathfrak{q}$.  (See \cite[pp. 288-289, 292]{h3}.) $\square$ \bigskip

Now consider the polynomial $G(x,j)$ from \cite[Section 2.2]{mor8}:
$$G(x,j) = (x^4-228x^3+494x^2+228x+1)^3-jx(1-11x-x^2)^5.$$
Recall from \cite[Prop. 5.5]{mor} that $\hat H_{5,l}(x)$ can be expressed as a product of factors $G(x,j)$ over the roots $j$ of $J_l(t)$, times certain factors corresponding to $j=0$ and $j=1728$:
\begin{equation}
\hat H_{5,l}(x) = c_{4,5}(x)^r (-c_{6,5}(x))^s x^{n_l}(1-11x-x^2)^{5n_l} J_l(j_5(x)),
\label{eq:20}
\end{equation}
with
\begin{align}
\label{eq:21} & j_5(x) = \frac{(x^4-228x^3+494x^2+228x+1)^3}{x(1-11x-x^2)^5},\\
\notag & c_{4,5}(x) = x^4-228x^3+494x^2+228x+1,\\
\notag & c_{6,5}(x) = -(x^2+1)(x^4+522x^3-10006x^2-522x+1).
\end{align}
Here $r, s, n_l$ and $J_l(t)$ have the same meaning as in the introduction.  We are using this polynomial instead of
$$F(x,j) = (x^4+12x^3+14x^2-12x+1)^3-jx^5(1-11x-x^2)$$
because of the identity (\ref{eq:13}).  The polynomial $G(x,j)$ will be especially useful in Section 6, since the roots of $G(x^5,j)$ are invariant under a group $G_{60} \cong A_5$ of linear fractional maps.  By (\ref{eq:20}) we can replace the rational function $j(x)$ in (\ref{eq:11}) by $j_5(x)$ (and multiply by the appropriate factor to clear denominators).  This yields exactly the same factors of the form $g(x)$ or $k(x)$, because of the fact that $j(\tau(x)) = j(\bar \tau(x)) = j_5(x)$: by (\ref{eq:7}) and (\ref{eq:9}) the polynomials $g(x)$ and $k(x)$ (with coefficients satisfying $r = \varepsilon^5(s-1)$) are essentially fixed by the map $(x \rightarrow \tau(x))$, while $(x \rightarrow \bar \tau(x))$ sends $k(x)$ to its companion polynomial $\bar k(x)$ (see the proof of Theorem 4.3(a)).  \medskip

We now prove the main result of this section. \bigskip

\noindent {\bf Theorem 5.4.} (a) {\it If $(X^2+a_iX+b_i)^2$ is a factor of $K_{5l}(X)$ in Theorem 2.1, whose roots are reductions mod $\mathfrak{p}$ $(\mathfrak{p} \mid \mathfrak{l}$ in $\Omega_f)$ of the $j$-invariants $\mathfrak{j}_i$ $(i = 1, 2)$, with $\mathfrak{j}_2 = \mathfrak{j}_1^{\sigma_{\wp_5}}$, then, for some $\sigma \in \textrm{Gal}(\Omega_f/K)$, the polynomial
$$g(x):= g_L(x)^\sigma \ \textrm{mod} \ \mathfrak{q}, \ \ (\mathfrak{p} \mid \mathfrak{q}, \ \mathfrak{q} \ \textrm{in} \ L),$$
is a factor of $G(x,\mathfrak{j}_1)G(x,\mathfrak{j}_2)$ mod $\mathfrak{q}$, and therefore divides $\hat H_{5,l}(x)$ over $\mathbb{F}_l$.}  \smallskip

(b) {\it If $l \equiv 2, 3 \ (\textrm{mod} \ 5)$, the polynomial $g(x) \in \mathbb{F}_l[x]$ is irreducible; while if $l \equiv 1, 4 \ (\textrm{mod} \ 5)$, $g(x)$ splits over $\mathbb{F}_l$ as a product of two polynomials of the form $k(x) = x^2 + rx +s$, with $r = \varepsilon^5(s-1)$ or $r = \bar \varepsilon^5(s-1)$, which are irreducible over $\mathbb{F}_l$.} \smallskip

(c) {\it If $g(x) \mid G(x,\mathfrak{k}_1)G(x,\mathfrak{k}_2)$ (mod $\mathfrak{q}$) for roots $\mathfrak{k}_1, \mathfrak{k}_2 = \mathfrak{k}_1^{\sigma_{\wp_5}}$ of $H_{-d}(X)$ and $d = 5l$ or $d = 20l$, then $\{\mathfrak{k}_1, \mathfrak{k}_2\} \equiv \{\mathfrak{j}_1, \mathfrak{j}_2\} \ (\textrm{mod} \ \mathfrak{p})$.  Thus, there is a unique factor of the form $(X^2+a_iX+b_i)^2$ in the congruence of Theorem 2.1 which gives rise to the factor $g(x)$ of $\hat H_{5,l}(x)$ over $\mathbb{F}_l$.} \medskip

\noindent {\it Proof.} Let $\mathfrak{j}_1=j(\tau_k)^\sigma$ and $\mathfrak{j}_2=j(\tau_k)^{\sigma_{\wp_5}\sigma}$ be two conjugate roots of $H_{-d}(X)$ ($d=5l$ or $20l$) over the field $L$ (for some $\sigma \in \textrm{Gal}(\Sigma_{\wp_5}\Omega_f/K)$).  Then $\rho_k^\sigma$ has the minimal polynomial $q(x, \lambda^\sigma)$ over $\Omega_f$ from the proof of Proposition 5.1.  The identity (\ref{eq:13}) implies that $\rho_k^\sigma=r^5(\tau_k)^\sigma$ is a root of $G(x,\mathfrak{j}_1) = 0$, so that $q(x,\lambda^\sigma) \mid G(x,\mathfrak{j}_1)$ in $\Omega_f[x]$.  Applying the automorphism $\sigma_{\wp_5}$ gives that $q(x,\lambda'^\sigma) \mid G(x,\mathfrak{j}_2)$, so that
$$g_L(x)^\sigma = q(x, \lambda^\sigma) q(x,\lambda'^\sigma) \mid G(x,\mathfrak{j}_1) G(x,\mathfrak{j}_2).$$
Let $\mathfrak{p}$ be a prime divisor of $\mathfrak{l}$ in $\Omega_f$, and $\mathfrak{q}$ the prime below $\mathfrak{p}$ in $L$.  If $\mathfrak{j}_1 \not \equiv \mathfrak{j}_2$ (mod $\mathfrak{p}$), and neither $\mathfrak{j}_i$ is congruent to $0$ or $1728$ (mod $\mathfrak{p}$), then $g_L(x)^\sigma \ \textrm{mod} \ \mathfrak{q} \in \mathbb{F}_l[x]$ is a quartic dividing $\hat H_{5,l}(x)$ (mod $\mathfrak{q}$), since $l \mid -d$ implies that the $j$-invariants $\mathfrak{j}_i$ reduce to supersingular $j$-invariants, which must be distinct roots of $J_l(t)$ over $\mathbb{F}_l$.  This proves (a).  \medskip

If $l \equiv 2,3$ (mod $5$), then the irreducible factors of $\hat H_{5,l}(x)$ (mod $l$) are either $x^2+1$ or quartic.  Now, $g_L(x)^\sigma \in L_\mathfrak{q}[x]$ is irreducible, and can only have a quadratic factor modulo $\mathfrak{q}$ if it is congruent to $(x^2+1)^2$.  But $\hat H_{5,l}(x)$ has no repeated factors, so $g_L(x)^\sigma$ must remain irreducible modulo $\mathfrak{q}$.  Hence, $g(x):= g_L(x)^\sigma \equiv x^4+ax^3+(11a+2)x^2-ax+1$ (mod $\mathfrak{q}$) is an irreducible factor of $\hat H_{5,l}(x)$ over $\mathbb{F}_l$. \medskip

In particular, this holds if $\mathfrak{j}_1$ and $\mathfrak{j}_2$ are two distinct roots mod $\mathfrak{p}$ of the factor $X^2+a_iX+b_i$ of $H_{-d}(X)$ ($d = 5l$ or $20l$) in Theorem 2.1.  This is because
\begin{equation}
\mathfrak{j}_2 \equiv \mathfrak{j}_1^l \equiv \mathfrak{j}_1^{\sigma_\mathfrak{l}} = \mathfrak{j}_1^{\sigma_{\wp_5}} \ (\textrm{mod} \ \mathfrak{p}),
\label{eq:22}
\end{equation}
so that
$$X^2+a_i X + b_i \equiv (X-\mathfrak{j}_1)(X-\mathfrak{j}_1^{\sigma_{\wp_5}}) \ (\textrm{mod} \ \mathfrak{p}).$$
Hence, $X^2+a_iX+b_i$ arises by reduction from roots $\mathfrak{j}_i$ which are conjugate over $L$.  It follows that the factor $(X^2+a_iX+b_i)^2$ contributes at least one irreducible factor of the form $g(x)$ to $\hat H_{5,l}(x)$ over $\mathbb{F}_l$, when $l \equiv 2, 3$ (mod $5$).  Note that different factors $X^2+a_iX+b_i$ will correspond to different automorphisms $\sigma \in \textrm{Gal}(\Omega_f/K)$ in the above discussion.  \medskip

On the other hand, if $l \equiv \pm 1$ (mod $5$), then the irreducible factors of $\hat H_{5,l}(x)$ are linear or quadratic.  Since $g_L(x)^\sigma$ mod $\mathfrak{q}$ divides $\hat H_{5,l}(x)$, then $g_L(x)^\sigma$ splits either as a product of two distinct quadratics or four distinct linear polynomials (mod $\mathfrak{q}$), by the proof of Proposition 4.2.  However, the latter cannot happen, by Hensel's Lemma, because $g(x)^\sigma$ is a product of irreducible quadratics over $L_\mathfrak{q}$ (Theorem 5.3(c)).  Hence, $g_L(x)^\sigma$ factors as a product of two irreducible quadratics (mod $\mathfrak{q}$).  Now setting $\lambda_1=\lambda^\sigma$ and $a = j_5^*(\tau_k)^\sigma = 22+\lambda_1+5^3 \lambda_1^{-1}$ gives
\begin{equation}
a^2-44a-16 = \frac{(\lambda_1^2-125)^2}{\lambda_1^2}.
\label{eq:23}
\end{equation}
By (\ref{eq:19}),
\begin{equation}
\left(\lambda_1-\frac{125}{\lambda_1}\right)^{\sigma_\mathfrak{l}}=\left(\lambda_1-\frac{125}{\lambda_1}\right)^{\sigma_\mathfrak{a}}=-\left(\lambda_1-\frac{125}{\lambda_1}\right).
\label{eq:24}
\end{equation}
Since $\sigma_\mathfrak{l}$ is the Frobenius automorphism for $\mathfrak{l}$ in $\Omega_f/K$, (\ref{eq:23}) and (\ref{eq:24}) imply that the quadratic residue symbol
\begin{equation*}
\left(\frac{-4\Theta_1}{\mathfrak{q}}\right) = \left(\frac{a^2-44a-16}{\mathfrak{q}}\right) = -1 \ \ \textrm{in} \ L.
\end{equation*}
The terms on either side of (\ref{eq:24}) are nonzero (mod $\mathfrak{q}$), because $a^2-44a-16$ is a factor of the discriminant of $g_L(x)^\sigma$; and we have shown above that $g_L(x)^\sigma$ has no multiple factor (mod $\mathfrak{q}$).  Hence, in the notation of Section 3, $-\Theta_1$ is not a square mod $\mathfrak{q}$ in $L$.  Now the proof of Proposition 4.2 shows that $g_L(x)^\sigma \equiv k_1(x) k_2(x)$ (mod $\mathfrak{q}$) factors into two polynomials of the form  $k_i(x) =x^2+r_ix+s_i$, where $r_i, s_i$ satisfy the conditions of that proposition. \medskip

Hence, when $l \equiv \pm 1$ (mod $5$), there are at least two distinct factors of $\hat H_{5,l}(x)$ of the form $k(x)$ arising from the factor $(X^2+a_iX+b_i)^2$ in Theorem 2.1. This completes the proof of (b). \medskip

The same arguments as in \cite[Lemma 4.2]{mor7} show that $H_{-5l}(X)$ and $H_{-20l}(X)$ are both squares (mod $l$), in the case $l \equiv 3$ (mod $4$), so that the factor $(X^2+a_iX+b_i)^2$ only divides one of these class equations (mod $l$).  Given this, the proof of (c) is straightforward.  If $\rho$ is a root of $g_L(x)^\sigma$ in $\Sigma_{\wp_5} \Omega_f$, then the hypothesis of (c) implies that $j_5(\rho) \equiv \mathfrak{k}_1$, say, modulo a prime divisor $\mathfrak{r}$ of $\mathfrak{p}$ in $\Sigma_{\wp_5} \Omega_f$.  But by (a), $j_5(\rho) \equiv \mathfrak{j}_i$ mod $\mathfrak{r}$, for $i = 1$ or $2$, so $\mathfrak{k}_1 \equiv \mathfrak{j}_i$ mod $\mathfrak{p}$.  The assertion follows from the relations $\mathfrak{j}_2 = \mathfrak{j}_1^{\sigma_{\wp_5}}$ and $\mathfrak{k}_2 = \mathfrak{k}_1^{\sigma_{\wp_5}}$.  $\square$

\section{Completion of the proof of Theorem 1.1}

It remains to show that, aside from the factors $g(x)$ and $k(x)$ in Theorem 5.4(b), there are no additional factors of the form $g(x)$ ($l \equiv 2, 3$ mod $5$) or $k(x)$ ($l \equiv 1, 4$ mod $5$) arising -- in the sense of Theorem 3.5 or equation (\ref{eq:11}) -- from a given factor $(X^2+a_iX+b_i)^2$ of $H_{-d}(X)$ in the congruence of Theorem 2.1, for $d=5l$ or $d=20l$.

\subsection{Sporadic quadratic factors for $d = 84, 96$.}

We do this first for the sporadic quadratic factors of $H_{-d}(X)$ mod $l$, which occur for $d = 84, 96$ and primes $l$ for which $\epsilon_d = 0$ in the proof of Theorem 2.1.  These are the quadratic factors of $H_{-d}(X)$ which divide $K_{5l}(X)$ only to the second power modulo $l$.  The factors of the form $g(x)$ or $k(x)$ which arise from these quadratic factors divide the polynomial
\begin{equation}
F_d(x)=x^{5\textsf{h}(-d)}(1-11x-x^2)^{\textsf{h}(-d)}H_{-d}(j(x)) \ \ \textrm{mod} \ l
\label{eq:25}
\end{equation}
from Theorem 3.5.  From \cite[Prop. 4.1]{mor2} we have the factorization
\begin{equation}
F_d(x) = Q_d(x) f_d(x), \ \ \textrm{deg}(Q_d(x)) = 4\textsf{h}(-d).
\label{eq:26}
\end{equation}
Considering the case $d=84$ in the proof of Theorem 2.1 and the non-bold primes $l > 379$ in Table 3 with $\epsilon_d = 0$, we must show that these factors, so far as they occur for $l \in \{389, 397, 401, 409 \}$ (see (\ref{eq:2}) and (\ref{eq:3})), each yield only one factor of the form $g(x)$ (or two of the form $k(x)$).  (The prime $383$ can be ignored for $d=84$ in Table 3, since it does not divide $N(Q(u,v))$ in equations (\ref{eq:2}) and (\ref{eq:3}).)  For example, the unique factors of the form $g(x)$ dividing $F_{84}(x)$ corresponding to these four primes are:
\begin{align*}
g_{389}(x) & = (x^2 + 286x + 379)(x^2 + 262x + 350), \ r \equiv 151^5(s-1) \ \textrm{mod} \ 389;\\
g_{397}(x) & = x^4 + 253x^3 + 6x^2 + 144x + 1;\\
g_{401}(x) & = (x^2 + 376x + 362)(x^2 + 205x + 329), \ r \equiv 111^5(s-1) \ \textrm{mod} \ 401;\\
g_{409}(x) & = (x^2 + 251x + 304)(x^2 + 240x + 74),  \ r \equiv 129^5(s-1) \ \textrm{mod} \ 409.
\end{align*}
Each of these polynomials divides the cofactor $f_{84}(x)$ of $Q_{84}(x)$ in $F_{84}(x)$, in the above notation.  In particular, this shows that only one of the quadratic factors of $H_{-84}(X)$ over $\mathbb{F}_l$ can divide $K_{5l}(X)$ for these primes, since we know from Theorem 5.4 that each distinct quadratic factor would yield a different $g(x)$ or pair of polynomials $k(x)$ dividing $\hat H_{5,l}(x)$.  (Note: the only primes in the set $\mathcal{S}$ which are listed for $d=84$ in Table 3 are $167, 227, 311$, and $H_{-84}(X)$ is a product of linear factors modulo each of them.) \medskip

A similar analysis applies to quadratic factors of $H_{96}(X)$ over $\mathbb{F}_l$, for the primes $l \in \{383, 397, 401, 421, 431, 449\}$ in Table 3, since $\epsilon_{96}=0$ for these primes.  (See equations (\ref{eq:4}) and (\ref{eq:5}).)  The primes $383, 431$ can be ignored, because $H_{-96}(X)$ splits completely for them, as it does for the primes $167, 239, 263, 359$; the latter are the primes in $\mathcal{S}$ listed in Table 3 for $d=96$.  We find the following unique factors of the form $g(x)$ of $F_{96}(x)$:
\begin{align*}
g_{397}(x) & =x^4 + 30x^3 + 332x^2 + 367x + 1;\\
g_{401}(x) & = (x^2 + 79x + 102)(x^2 + 184x + 287), \ r \equiv 289^5(s-1) \ \textrm{mod} \ 401;\\
g_{421}(x) & = (x^2 + 316x + 351)(x^2 + 209x + 6), \ r \equiv 110^5(s-1) \ \textrm{mod} \ 421;\\
g_{449}(x) & = (x^2 + 92x + 307)(x^2 + 437x + 332),  \ r \equiv 165^5(s-1) \ \textrm{mod} \ 449.
\end{align*}
This calculation yields: \medskip

\noindent {\bf Lemma 6.1.} {\it If $d = 84$ or $96$, and $X^2+a_iX+b_i$ is a sporadic factor of $H_{-d}(X)$ mod $l$ in the congruence of Theorem 2.1, for some prime $l >379$ for which $\epsilon_d = 0$, then there is a unique factor of the form $g(x)$ dividing $\hat H_{5,l}(x)$ over $\mathbb{F}_l$ which also divides $F_d(x)$, i.e. which corresponds to the factor $(X^2+a_iX+b_i)^2$.  There are no sporadic quadratic factors of $H_{-d}(X)$ for primes $l$ lying in the set $\mathcal{S}$.} \medskip
 
Recall that the non-sporadic factors of $H_{-d}(X)$ (for $d= 84, 96$) have already been discussed in Section 4.  See especially Propositions 4.1 and 4.3.

\subsection{General factors $X^2+a_iX+b_i$ of $K_{5l}(X)$ mod $l$.}

Next, we state the details concerning the icosahedral group that we will need from \cite[\S 2.2]{mor8} for the rest of the argument.  The letter $\zeta$ denotes a fixed primitive $5$-th root of unity.  \medskip

\noindent {\bf Lemma 6.2.} (a) {\it The linear fractional maps
\begin{align*}
S(x) = \zeta x,  \ \ & T(x) = \frac{-(1+\sqrt{5})x+2}{2x+1+\sqrt{5}},\\
U(x) = \frac{-1}{x}, \   & A(x) =\zeta^3 \frac{(1+\zeta)x+1}{x-1-\zeta^4}
\end{align*}
satisfy the relations
\begin{equation*}
S^5(x) = x, \ \ T^2(x) = x, \ \ U^2(x) = x, \ \ A^3(x) = x;
\end{equation*}
and
\begin{align*}
A &= STS^{-2}, \ \ \ \ \ \ A^\sigma = A^{-1}U,\\
ATA^{-1} &= U, \ \ \ AUA^{-1} = TU = UT = T_2;
\end{align*}
here $\sigma = (\zeta \rightarrow \zeta^2)$ and $A^\sigma$ is the result of applying $\sigma$ to the coefficients of $A$.} \smallskip

\noindent (b) {\it The group $G_{60} = \langle S, T \rangle$ is Fricke's normal form of the icosahedral group \cite[II, pp. 41-43]{fr}.} \smallskip

\noindent (c) {\it The group $H = \langle T, U \rangle = \{1, T, U, T_2\}$ is isomorphic to the Klein $4$-group and $\langle H, A \rangle \cong A_4$, the alternating group on $4$ letters.  Further, the elements $S^i A^k$ are representatives of the left cosets of $H$ in $G_{60}$.} \smallskip

\noindent (d) {\it The group $G_{10}=\langle S, U \rangle$ generated by $S$ and $U$ has order $10$.  Representatives of the right cosets of $G_{10}$ in $G_{60}$ are the elements $T^i A^k$, for $i = 0, 1$ and $k = 0, 1, 2$.}
\medskip

\noindent {\it Proof.} Parts (a), (b), (c) are discussed in \cite[\S 2.2]{mor8} and are left to the reader.  To prove (d), we just note the relations
\begin{align*}
AU & = UTA,\\
AT & = UA,\\
AT_2 & = AUT = TA.
\end{align*}
By (c), every element of $G_{60}$ lies in one of the cosets $S^i A^k H$.  Using this, the above relations imply that every element of $G_{60}$ lies in one of the cosets $G_{10}T^i A^k$, and the assertion of (d) follows. $\square$ \medskip

We turn now to the non-sporadic factors $(X^2+a_i X + b_i)^2$ of $K_{5l}(X)$ in Theorem 2.1.  We shall use the fact that a factor of $\hat H_{5,l}(x)$ of the form $k(x) = x^2 +rx+s$, satisfying $r=\varepsilon^5(s-1)$, has (nonzero) roots $\alpha^5$ and $\beta^5$ satisfying the equation
\begin{equation}
\alpha^5+\beta^5 = \varepsilon^5(1-\alpha^5 \beta^5).
\label{eq:27}
\end{equation}
(Theorem 5.8 of \cite{mor} shows that $\alpha, \beta$ lie in $\mathbb{F}_{l^2}$, resp. $\mathbb{F}_{l^4}$, according as $l \equiv \pm 1$ or $l \equiv \pm 2$ mod $5$.)  Furthermore, the quantity $\beta^5$ is determined by $\alpha^5$, and vice versa, since
$$\beta^5 = \frac{-\alpha^5+\varepsilon^5}{\varepsilon^5 \alpha^5 +1} = \tau(\alpha^5).$$
Since there is at least one such factor $k(x)$ arising from $(X^2+a_i X + b_i)^2$ (which may divide an irreducible $g(x)$), we can choose a fixed $\alpha$ satisfying (\ref{eq:27}), which is unique up to multiplication by a power of $\zeta$, for which $\rho=\alpha^5$ is a root of $k(x)$ or $g(x) = k_1(x) k_2(x)$.  Since, moreover, $(\alpha, \beta)$ and $(-1/\alpha,-1/\beta)$ are also solutions of this equation, it is clear that for $M$ in the group $G_{10}=\langle S, U \rangle$, the pairs $(M(\alpha), M(\beta))$ are solutions of (\ref{eq:27}) for which $M(\alpha)^5 = \alpha^5$ or $-1/\alpha^5$, the latter quantity being a root of the companion polynomial
$$\bar k(x) = \frac{1}{s} x^2 k(-1/x) = x^2-\frac{r}{s}x+\frac{1}{s}$$
and therefore a root of $g(x) = k(x) \bar k(x)$, in any case.  Thus, ``conjugating'' $(\alpha, \beta)$ by an element in the group $G_{10}$ gives a solution of (\ref{eq:27}) corresponding to the same factor $g(x)$ or pair of factors $k(x), \bar k(x)$. \medskip

Further, by Theorem 5.4(a) and (\ref{eq:22}), $k(x)$ divides a product
$$G(x,\mathfrak{j}_1)G(x,\mathfrak{j}_2) \equiv G(x,\mathfrak{j}_1)G(x,\mathfrak{j}_1^l) \ \textrm{mod} \ \mathfrak{p},$$
so $\alpha^5$ must be a root of $G(x, j) = 0$ or $G(x, j^l) = 0$, for the reduced $j$-invariant $j \equiv \mathfrak{j}_1$ (mod $\mathfrak{p}$), say.    
We introduce the following terminology: we say that
$$\alpha^5 \ belongs \ to \ the \ invariant \ j \ if \ j = j_5(\alpha^5), i.e. \ if \ G(\alpha^5,j) = 0 \ in \ \overline{\mathbb{F}}_l;$$
and $j_5(x)$ is the rational function in (\ref{eq:21}).  Clearly, a root $\alpha^5$ of a factor $g(x)$ or $k(x)$ belongs to a unique $j$-invariant $j$.  If $l \equiv 1, 4$ (mod $5$) and $\alpha^5, \beta^5$ are the roots of the quadratic $k(x) \in \mathbb{F}_l[x]$, then $\beta^5 = \alpha^{5l}$.  Hence, if $\alpha^5$ belongs to $j$, $\beta^5$ belongs to $j^l$.  If $l \equiv 2, 3$ (mod $5$), and $\alpha^5, \beta^5$ are roots of the quartic $g(x) \in \mathbb{F}_l[x]$ satisfying (\ref{eq:27}), then $\beta^5 = \tau(\alpha^5)$, which equals $\alpha^{5l}$ or $\alpha^{5l^3} = \frac{-1}{\alpha^{5l}} = U(\alpha^{5l})$, as in the discussion leading up to Theorem 3.4.  Since the map $x \rightarrow U(x)$ fixes the rational function $j_5(x)$, we also see in this case that if $\alpha^5$ belongs to $j$, then $\beta^5$ belongs to $j^l$.  \medskip
  
Thus, a second factor of the form $g(x)$ or $k(x)$ (distinct from $\bar k(x)$) dividing $G(x,j) G(x,j^l)$ would yield a root $\alpha'^5$ of $G(x,j) = 0$ (w.l.o.g.) and a solution $(\alpha',\beta')$ of (\ref{eq:27}) or of the equation
\begin{equation*}
\alpha'^5+\beta'^5 = \bar \varepsilon^5(1-\alpha'^5 \beta'^5).
\end{equation*}
Now, the $60$ roots of $G(x^5,j) = 0$ have the form $M(\alpha)$, for a single root $\alpha$ and $M \in G_{60}$ (see \cite[p. 243]{mor8}), and this is true for any $j$.  Since the pairs $\{\alpha^5, \alpha'^5\}$ and $\{\beta^5, \beta'^5\}$ belong to the invariants $j$ and $j^l$, respectively, we see that there are $M_1, M_2 \in G_{60}$ for which:
\begin{equation}
\alpha'=M_1(\alpha) \ \textrm{is not in the orbit} \ G_{10}\alpha \ \textrm{and} \ \beta'=M_2(\beta) \notin G_{10}\beta.
\label{eq:28}
\end{equation}
Note that $M_1(\alpha) \in G_{10}\alpha$ if and only if $M_2(\beta) \in G_{10}\beta$, since $\alpha^5$ and $\beta^5$ determine each other.  
\medskip

To prove that {\it only} one factor of the form $g(x)$, or two factors of the form $k(x)$, correspond to the factor $(X^2+a_iX+b_i)^2=(X-j)^2(X-j^l)^2$, we consider the equations
\begin{align}
\label{eq:29} M_1(\alpha)^5+M_2(\beta)^5 & =\varepsilon^5(1-M_1(\alpha)^5 M_2(\beta)^5),\\
\label{eq:30} M_1(\alpha)^5+M_2(\beta)^5& = \bar \varepsilon^5(1-M_1(\alpha)^5 M_2(\beta)^5).
\end{align}
By (\ref{eq:28}) it suffices to consider these equations for pairs of elements $(M_1, M_2)$ with $M_i \notin G_{10}$.  Hence we can take both $M_i$ to be of the form $M = T^i A^k$, where $i, k$ are not both $0$ (Lemma 6.2(d)).  Considering (\ref{eq:27}) and (\ref{eq:29}), we first compute the resultants
\begin{align*}
R_{M_1,M_2} =& \ \textrm{Res}_y(x^5+y^5-\varepsilon^5(1-x^5 y^5),\\
& (c_1x+d_1)^5(c_2y+d_2)^5(M_1(x)^5+M_2(y)^5-\varepsilon^5(1-M_1(x)^5 M_2(y)^5))),
\end{align*}
where $M_i(x) = \frac{a_ix+b_i}{c_ix+d_i}$.  For example,
\begin{align*}
R_{T,T} &= 5^{15} x (x^2+x-1)(x^2+1)(x^4 - x^3 + x^2 + x + 1)\\
& \ \ \times (x^4 - 2x^3 + 2x + 1) (x^4 + x^3 + 3x^2 - x + 1) \\
& \ \ \times (x^8 + 4x^7 + 10x^6 + 8x^5 + 12x^4 - 8x^3 + 10x^2 - 4x + 1) \\
& \ \ \times (x^8 + 7x^7 + 15x^6 + 15x^5 + 16x^4 - 15x^3 + 15x^2 - 7x + 1)\\
& \ \ \times (x^{16} + 2x^{15} - 4x^{14} - 12x^{13} + 25x^{12} - 18x^{11} + 68x^{10} - 112x^9\\
& \ \ \ \ + 13x^8 + 112x^7 + 68x^6 + 18x^5 + 25x^4 + 12x^3 - 4x^2 - 2x + 1)\\
& = \ 5^{15} x (x^2+x-1) p_4(x) p_{11}(x) p_{16}(x) p_{19}(x)p_{64}(x) p_{99}(x) p_{84}(x),
\end{align*}
where $p_d(x)$ is the polynomial defined in \cite{mor2} (see Tables 1 and 2 in that paper) and $p_4(x)=x^2+1$.  By the results of \cite[pp. 1193-1195]{mor2}, the roots of $p_d(x)$ are solutions of (\ref{eq:27}) {\it in characteristic zero}, and $p_d(x)$ divides $Q_d(x^5) = p_d(x)q_d(x)$; so by (\ref{eq:26}) $p_d(x)$ also divides
$$F_d(x^5)=x^{25\textsf{h}(-d)}(1-11x^5-x^{10})^{\textsf{h}(-d)}H_{-d}(j(x^5)).$$
The factors $x(x^2+x-1)$ divide $x(x^{10}+11x^5-1)$ and can be ignored, since their roots correspond to singular curves $E_5(b)$.  This calculation implies that the solution $(\alpha,\beta)$ must arise from one of the factors $H_{-d}(X)^{4\epsilon_d}$, for $d \in \mathfrak{T}$.  If $X^2+a_iX+b_i$ is not a sporadic quadratic, as discussed in section 6.1, and therefore does not divide any of the class equations $H_{-d}(X)$ (mod $l$), for $d \in \mathfrak{T}$, this shows that no solution of (\ref{eq:27}) of the form $(\alpha', \beta') = (T(\alpha),T(\beta))$ can arise from the same factor $(X^2+a_iX+b_i)^2$ that $(\alpha, \beta)$ does.  \medskip

For the map $A(x)= \zeta^3 \frac{(1+\zeta)x+1}{x-1-\zeta^4}$, the resultant $R_{A,A}$ above is a product of polynomials in $\mathbb{Q}(\zeta)[x]$ whose norm to $\mathbb{Q}$ is:
\begin{align*}
N_\mathbb{Q}(R_{A,A})= & \ 5^{60}x^4 (x^4 - 3x^3 + 4x^2 - 2x + 1)(x^4 + 2x^3 + 4x^2 + 3x + 1)\\
& \times q_4(x) q_{11}(x) q_{16}(x) q_{19}(x) q_{64}(x) q_{84}(x) q_{99}(x),
\end{align*}
where $Q_d(x^5) = p_d(x)q_d(x)$, as in \cite[p. 1195]{mor2}.  We have the relation $q_d(x) = \prod_{i=1}^4{p_d(\zeta^i x)}$, with $\zeta$ as in Lemma 6.2.  The factors $q_d(x)$ also divide the polynomial $F_d(x^5)=Q_d(x^5)f_d(x^5)$, so the same arguments apply as above.  (The two quartic factors divide $x^{10}+11x^5-1$ and can be ignored.)  \medskip

Since $N_\mathbb{Q}(R_{A^2,A^2}) = N_\mathbb{Q}(R_{TA,TA})=N_\mathbb{Q}(R_{TA^2,TA^2})=N_\mathbb{Q}(R_{A,A})$, we don't get any new factors from these pairs.  We still have to check the pairs $(M_1,M_2)$ with $M_1 \neq M_2$.  It suffices to check $10$ pairs with $M_1 \neq M_2$, since $(M_1,M_2)$ and $(M_2,M_1)$ yield the same solutions.  We find that the resultants
$$R_{T,A}=R_{T,TA}=R_{T,A^2}$$
coincide with $R_{T,T}$, so these pairs don't give anything new.  Also,
\begin{align*}
R_{T,TA^2}&=5^{15}\varepsilon^{25}(x^2+1)(x^8 - 2x^7 + x^6 - 4x^5 + 3x^4 + 4x^3 + x^2 + 2x + 1)\\
& \ \ \times (x^8 + x^6 - 6x^5 + 9x^4 + 6x^3 + x^2 + 1)\\
& \ \ \times (x^8 + x^7 + x^6 - 7x^5 + 12x^4 + 7x^3 + x^2 - x + 1)\\
& \ \ \times (x^8 + 4x^7 - x^6 - 14x^5 + 23x^4 + 14x^3 - x^2 - 4x + 1)\\
& \ \ \times (x^{16} + 4x^{15} + 29x^{12} - 24x^{11} + 86x^{10} - 32x^9 + 105x^8 + 32x^7\\
& \ \ \ \ \ + 86x^6 + 24x^5 + 29x^4 - 4x + 1)\\
& = 5^{15}\varepsilon^{25} p_4(x) p_{24}(x) p_{36}(x) p_{51}(x) p_{91}(x) p_{96}(x),\\
N_\mathbb{Q}(R_{A,TA^2})&=5^{60}q_4(x)q_{24}(x)q_{36}(x)q_{51}(x)q_{91}(x)q_{96}(x);
\end{align*}
so that these two products account for the remaining integers $d$ in the set $\mathfrak{T}$.  Furthermore, the resultant norms
$$N_\mathbb{Q}(R_{A,TA})=N_\mathbb{Q}(R_{A,A^2})=N_\mathbb{Q}(R_{A^2,TA})$$
coincide with $N_\mathbb{Q}(R_{A,A})$ above, while the resultant norms
$$N_\mathbb{Q}(R_{A^2,TA^2})=N_\mathbb{Q}(R_{TA,TA^2})$$
coincide with $N_\mathbb{Q}(R_{A,TA^2})$.  This accounts for all $10$ pairs.  \medskip

We must also account for possible simultaneous solutions of (\ref{eq:27}) and (\ref{eq:30}) and therefore must also consider the resultants
\begin{align*}
\bar R_{M_1,M_2} =& \textrm{Res}_y(x^5+y^5-\varepsilon^5(1-x^5 y^5),\\
& (c_1x+d_1)^5(c_2y+d_2)^5(M_1(x)^5+M_2(y)^5- \bar \varepsilon^5(1-M_1(x)^5 M_2(y)^5))).
\end{align*}
We check that this yields only the same solutions as before:
\begin{align*}
\bar R_{T,T}= & -5^{15} \bar \varepsilon^{25} p_4(x) p_{24}(x) p_{36}(x) p_{51}(x) p_{91}(x) p_{96}(x),\\
N_\mathbb{Q}(\bar R_{A,A})= & \ 5^{60} q_4(x) q_{24}(x) q_{36}(x) q_{51}(x) q_{91}(x) q_{96}(x),\\
N_\mathbb{Q}(\bar R_{TA,TA}) = & \ N_\mathbb{Q}(\bar R_{A^2,A^2}) = N_\mathbb{Q}(\bar R_{TA^2,TA^2})\\
=& \ N_\mathbb{Q}(\bar R_{A,A});\\
\bar R_{T,A} = & \ \bar R_{T,TA} =  \bar R_{T,A^2} = \bar R_{T,T}, \ \bar R_{T,TA^2} = -R_{T,T};\\
N_\mathbb{Q}(\bar R_{A,TA}) = & N_\mathbb{Q}(\bar R_{A,A^2}) = N_\mathbb{Q}(\bar R_{A^2,TA}) = N_\mathbb{Q}(\bar R_{A,A});\\
N_\mathbb{Q}(\bar R_{A,TA^2}) = & N_\mathbb{Q}(\bar R_{A^2,TA^2}) = N_\mathbb{Q}(\bar R_{TA,TA^2}) = N_\mathbb{Q}(R_{A,A}).
\end{align*}

When $l \equiv \pm 1$ (mod $5$) we would need to do the same calculation for the polynomials $k(x) = x^2+rx+s$ with $r = \bar \varepsilon^5(s-1)$, and work with the conjugate equation
\begin{equation}
\alpha^5+\beta^5 = \bar \varepsilon^5(1-\alpha^5 \beta^5)
\label{eq:31}
\end{equation}
in place of (\ref{eq:27}).  But these calculations follow from what we have already computed, since we can just apply the automorphism $\sigma = (\zeta \rightarrow \zeta^2)$ to the maps in $G_{60}$, and this map switches (\ref{eq:27}) and (\ref{eq:31}).  This sends the group $G_{10}$ to itself, and replaces $T$ by $T^\sigma = T_2 = UT$ and $A$ by $A^\sigma = A^{-1}U=TA^2$.  Since $T^\sigma A^\sigma = T_2 A^{-1}U = T_2 A^2 U = T_2 T A^2 = UA^2$, and $T_2A^k$ = $UTA^k$, we obtain exactly the same polynomials on taking norms as before. \medskip

Combining this discussion with Lemma 6.1, we have the following. \medskip

\noindent {\bf Proposition 6.3} {\it Each of the factors $(X^2+a_iX+b_i)^2$ in Theorem 2.1 contributes exactly one irreducible factor of the form $g(x)$ ($l \equiv 2, 3$ mod $5$) or exactly two irreducible factors of the form $k(x)$ ($l \equiv 1, 4$ mod $5$) to the factorization of $\hat H_{5,l}(x)$ over $\mathbb{F}_l$, for $l \in \mathcal{S}$ or $l >379$.}  \medskip

By the arguments in the final paragraphs of Sections 4.1 and 4.2, this proves Theorem 1.1 for the primes $l \in \mathcal{S} \cup \{l: l>379\}$. \medskip

The counts of quartic and quadratic factors for the primes satisfying $7 \le l \le 379$ and $l \not \in \mathcal{S}$ are given in Tables 6-9.  The number $N$ of such factors was counted by hand, after computing $\hat H_{5,l}(x)$ on Maple.  In each case, $N$ agrees with the formula of Theorem 1.1.  This shows that Theorem 1.1 holds for all primes $l >5$.

\begin{table}
  \centering 
  \caption{Number $N$ of factors $g(x)$ or $k(x)$ modulo $l \equiv 1$ mod $12$.}\label{tab:5}
  
\noindent \begin{tabular}{|c|l|c|c|}
\hline
   &  &  &  \\
$l$	&   $N$  & $h(-5l)$ & formula of Thm. 1.1\\
\hline
 13 & 2 & 8 & 2\\
 37 & 4 & 16 & 4\\
 61  & 8  &  16 & 8 \\
73  &  5  & 20 & 5\\
 97 & 5 & 20 & 5 \\
 109 & 16 & 32 & 16\\
  157 & 4 & 16 & 4\\ 
  181 & 12 & 24 &  12\\
  229 & 12 & 24 & 12\\
  241 & 20 & 40 & 20\\
  277 & 12 & 48 & 12\\
  313 & 7 & 28 & 7\\
  349 & 20 & 40 & 20\\
  373 & 12 & 48 & 12\\
    \hline
\end{tabular}
\end{table}

\begin{table}
  \centering 
  \caption{Number $N$ of factors $g(x)$ or $k(x)$ modulo $l \equiv 5$ mod $12$.}\label{tab:6}
  
\noindent \begin{tabular}{|c|l|c|c|}
\hline
   &  &  &  \\
$l$	&   $N$  & $h(-5l)$ & formula of Thm. 1.1\\
\hline
 17 & 1 & 4 & 1\\
 29 & 4 & 8 & 4\\
 41  & 4  &  8 & 4 \\
53  &  2  & 8 & 2\\
 89 & 4 & 8 & 4 \\
 113 & 3 & 12 & 3\\
  137 & 3 & 12 & 3\\ 
  149 & 8 & 16 &  8\\
  197& 6 & 24 & 6\\
  233 & 5 & 20 & 5\\
  257 & 3 & 12 & 3\\
  281 & 12 & 24 & 12\\
  317 & 6 & 24 & 6\\
  353 & 5 & 20 & 5\\
    \hline
\end{tabular}
\end{table}

\begin{table}
  \centering 
  \caption{Number $N$ of factors $g(x)$ or $k(x)$ modulo $l \equiv 7$ mod $12$.}\label{tab:7}
  
\noindent \begin{tabular}{|c|l|c|c|c|}
\hline
   &  &  &  &  \\
$l$	&   $N$  & $h(-5l)$ & $l$ mod $8$ & formula of Thm. 1.1\\
\hline
 7 & 1 & 2 & 7 & 1\\
 19 & 5 & 8 & 3 & 5\\
 31  & 7  &  4 & 7 & 7 \\
43  &  6  & 14 & 3 & 6\\
 67 & 8 & 18 & 3 & 8 \\
 79 & 13 & 8 & 7 & 13\\
  127 & 9 & 10 & 7 & 9\\ 
  139 & 21 & 24 & 3 & 21\\
  151& 23 & 12 & 7 &23\\
  163 & 14 & 30 & 3 &14\\
  211 & 35 & 36 & 3 &35\\
  283 & 16 & 34 & 3 &16\\
  307 & 18 & 38 & 3 &18\\
  331 & 43 & 44 & 3 &43\\
  379 & 45 & 48 & 3 &45\\
  \hline
\end{tabular}
\end{table}

\begin{table}
  \centering 
  \caption{Number $N$ of factors $g(x)$ or $k(x)$ modulo $l \equiv 11$ mod $12$.}\label{tab:8}
  
\noindent \begin{tabular}{|c|l|c|c|c|}
\hline
   &  &  &  &  \\
$l$	&   $N$  & $h(-5l)$ & $l$ mod $8$ & formula of Thm. 1.1\\
\hline
 11 & 3 & 4 & 3 &3\\
 23 & 1 & 2 & 7 &1\\
47  &  1  &  2 & 7 &1 \\
59  &  5  & 8 & 3  &5\\
 71 & 7 & 4 & 7 & 7 \\
 83 & 4 & 10 & 3 &4\\
  131 & 11 & 12 & 3 & 11\\ 
  \hline
\end{tabular}
\end{table}

\section{The degree of $ss_p^{(5*)}(X)$}

\begin{table}
  \centering 
  \caption{Degree and number of linear factors of $ss_p^{(5*)}(X)$.}\label{tab:9}
  
\noindent \begin{tabular}{|c|c|c||c|c|c|}
\hline
   &  &  &  &  & \\
$p$	&   $\textrm{deg}(ss_p^{(5*)}(X))$  & $L^{(5*)}(p)$ & $p$	&   $\textrm{deg}(ss_p^{(5*)}(X))$  & $L^{(5*)}(p)$\\
\hline
 7 & 2 & 2 & 47 & 12 & 2 \\
11 & 4 & 4 & 53 & 14 & 2 \\
13  &  4  &  2  & 59 & 16 & 10 \\
17  &  5  & 1 & 61 & 15 & 7 \\
 19 & 6 & 6 & 67 & 17 & 9 \\
 23 & 6 & 2 & 71 & 19 & 11 \\
 29 & 7 & 5 & 73 & 19 & 5 \\ 
 31 & 9 & 7 & 79 & 21 & 13 \\
 37 & 10 & 4 & 83 & 21 & 5 \\
 41 & 10 & 6 &  89 & 22 & 8 \\
 43 & 11 & 7 & 97 & 25 & 5 \\
  \hline
\end{tabular}
\end{table}

The aim of this section is to prove Nakaya's conjecture \cite[Conjecture 6]{na} for $N=5$, namely, that
$$\textrm{deg}(ss_p^{(5*)}(X)) = \frac{1}{4}\left(p-\left(\frac{-1}{p}\right)\right)+\frac{1}{2}\left(1-\left(\frac{-5}{p}\right)\right), \ \ p > 5.$$
Recall that $ss_p^{(5*)}(X)$ is determined as follows.  Define the polynomial $R_5(X,Y)$ by
\begin{align*}
R_5(X,Y) = &X^2-X(Y^5-80Y^4+1890Y^3-12600Y^2+7776Y+3456)\\
&+(Y^2+216Y+144)^3.
\end{align*}
Then
$$ss_p^{(5*)}(X) = \prod_{j_5^*}{(X-j_5^*)} \in \mathbb{F}_p[X],$$
where $j_5^*$ runs over the distinct roots of $R_5(j,j_5^*) \equiv 0$ in $\overline{\mathbb{F}}_p$ for those values of $j$ which are supersingular in characteristic $p$, i.e., the roots of $ss_p(X)=0$. \medskip

We use the parametrization
$$(j,j_5^*)=(X,Y)=\left(-\frac{(z^2+12z+16)^3}{z+11}, -\frac{z^2+4}{z+11}\right)$$
of $R_5(X,Y)=0$ given in \cite{mor8} (see the proof of Theorem 6.1).  Note first that
$$\textrm{disc}_z((z^2+12z+16)^3+j(z+11))=3125j^4(j - 1728)^2,$$
so that there are exactly $6$ values of $z$ for every supersingular value of $j$, except for $j=0$ and $j=1728$.  For these values
\begin{align*}
j & =0 \ \ \textrm{iff} \ \ z^2+12z+16=0,\\
j & = 1728 \ \ \textrm{iff} \ \ (z^2 + 4)(z^2 + 18z + 76)=0;
\end{align*}
where the two quadratics have discriminants $2^4 \cdot 5$ and $-2^{14} \cdot 3^8 \cdot 5^3$, respectively.  Thus exactly $2$, respectively $4$, values of $z$ correspond to $j=0$ and $j=1728$ in characteristic $p > 5$.  Since the number of supersingular $j$-invariants in characteristic $p$ is given by
$$n_p+r_p+s_p=\frac{p-e_p}{12}+\frac{1}{2}\left(1-\left(\frac{-3}{p}\right)\right)+\frac{1}{2}\left(1-\left(\frac{-1}{p}\right)\right),$$
where $e_p \in \{1,5,7,11\}$ and $p \equiv e_p$ (mod $12$), there are exactly
$$\frac{p-e_p}{2}+\left(1-\left(\frac{-3}{p}\right)\right)+2\left(1-\left(\frac{-1}{p}\right)\right)$$
values of the parameter $z$ altogether.  Using that $e_p=6-2\left(\frac{-3}{p}\right)-3\left(\frac{-1}{p}\right)$, we find
$$\frac{1}{2}\left(p-\left(\frac{-1}{p}\right)\right)$$
values of $z$ which correspond to supersingular $j$-invariants in $\mathbb{F}_{p^2}$. \medskip

Now note that
$$\textrm{disc}_z(z^2+4+t(z+11))=t^2-44t-16, \ \ \textrm{disc}(t^2-44t-16)=2^4 \cdot 5^3.$$
It follows that exactly two values of $z$ correspond to a single value of $j_5^*$, except when $t=j_5^*$ is a root of the last quadratic.  Since
$$\text{Res}_t(z^2+4+t(z+11),t^2-44t-16) = (z^2+22z-4)^2,$$
the roots of $z^2+22z-4$ correspond $1-1$ to the roots $j_5^*$ of $t^2-44t-16=0$.  On the other hand, the corresponding values of $j$ are supersingular exactly when they are roots of
\begin{align*}
Res_z((z^2 + 12z + 16)^3 &+ j(z + 11), z^2 + 22z - 4)\\
& =-5^3(j^2 - 1264000j - 681472000)\\
& = -5^3 H_{-20}(j).
\end{align*}
The roots of $H_{-20}(x)$ are supersingular in $\mathbb{F}_{p^2}$ if and only if $\left(\frac{-5}{p}\right) = -1$.  Therefore, altogether we have
$$\frac{1}{4} \left(p-\left(\frac{-1}{p}\right)\right)$$
distinct values of $j_5^*$ corresponding to supersingular $j$-invariants, when $\left(\frac{-5}{p}\right)=+1$, and
$$\frac{1}{2} \left(\frac{1}{2}\left(p-\left(\frac{-1}{p}\right)\right)-2\right)+2=\frac{1}{4} \left(p-\left(\frac{-1}{p}\right)\right)+1$$
distinct values of $j_5^*$ corresponding to supersingular $j$-invariants, when $\left(\frac{-5}{p}\right)=-1$.  Hence, the degree of $ss_p^{(5*)}(X)$ is given by
$$\frac{1}{4}\left(p-\left(\frac{-1}{p}\right)\right)+\frac{1}{2}\left(1-\left(\frac{-5}{p}\right) \right).$$
This proves Nakaya's Conjecture 6 for $N=5$. \medskip

\noindent {\bf Acknowledgement.} I am grateful to an anonymous referee for his careful reading of the manuscript and many useful suggestions, which have led to significant improvements in the exposition in several places.

\medskip

\noindent Dept. of Mathematical Sciences, LD 270

\noindent Indiana University - Purdue University at Indianapolis (IUPUI)

\noindent Indianapolis, IN 46202

\noindent {\it e-mail: pmorton@iupui.edu}

\end{document}